\documentclass[12pt]{amsart}

\usepackage{epsf}
\usepackage{amsmath,amssymb}
\usepackage{colordvi}
\usepackage{color}
\usepackage{graphics}

\newcommand{\map}[1]{\xrightarrow{#1}}

\theoremstyle{plain}

\newtheorem{thm}{Theorem}
\numberwithin{thm}{section}
\newtheorem{prop}[thm]{Proposition}
\newtheorem{lem}[thm]{Lemma}
\newtheorem{cor}[thm]{Corollary}

\theoremstyle{definition}
\newtheorem{defi}[thm]{Definition}
\newtheorem{rem}[thm]{Remark} %use \begin{rem}

\newcounter{FNC}[page]
\def\newfootnote#1{{\addtocounter{FNC}{2}$^\fnsymbol{FNC}$%
     \let\thefootnote\relax\footnotetext{$^\fnsymbol{FNC}$#1}}}
\newcommand{\sumsub}[1]{\sum_{\substack{#1}}} % sums with multline subscripts

\newcommand{\ten}{\mbox{\raisebox{1pt}{${\scriptstyle \otimes}$}}}
\newcommand{\iten}{{\backslash}}       % for internal tensors
\newcommand{\st}{\mathrm{st}}          % for standarization
\newcommand{\id}{\mathit{id}}          % for identity

\newcommand{\Inv}{\mathrm{Inv}}

\newcommand{\Des}{\mathrm{Des}}
\newcommand{\GDes}{\mathrm{GDes}}             %for (proper) global descents
\newcommand{\gr}{\rm{gr}}

\newcommand{\setF}{F}%{\mathsf{F}}
\newcommand{\setM}{M}%{\mathsf{M}}
\newcommand{\setR}{\mathsf{R}}
\newcommand{\setS}{\mathsf{S}}
\newcommand{\setT}{\mathsf{T}}

\newcommand{\Q}{\mathbb{Q}}

\newcommand{\Sym}{\mathit{Sym}}
\newcommand{\QSym}{\mathcal{Q}\Sym}
\newcommand{\SSym}{\mathfrak{S}\Sym}
\newcommand{\YSym}{\mathcal{Y}\Sym}
\newcommand{\NSym}{\mathcal{N}\!\Sym}
\newcommand{\LR}{\mathit{LR}}
\newcommand{\CK}{\mathit{CK}}
\newcommand{\GL}{\mathit{GL}}
\newcommand{\NCK}{\mathit{NCK}}

\newcommand{\frakS}{\mathfrak{S}}
\newcommand{\Sh}[1]{\mathfrak{S}^{#1}}

\newcommand{\calA}{\mathcal{A}}
\newcommand{\calC}{\mathcal{C}}
\newcommand{\calD}{\mathcal{D}}

\newcommand{\calL}{\mathcal{L}}
\newcommand{\calF}{F}%{\mathcal{F}}
\newcommand{\calM}{M}%{\mathcal{M}}
\newcommand{\calQ}{\mathcal{Q}}

\newcommand{\calY}{\mathcal{Y}}
\newcommand{\calZ}{\mathcal{Z}}

\newcommand{\barBell}{\begin{picture}(5.5,12)(-2.5,2)\thicklines\put(0,0){\circle*{2.5}}
           \put(0,0){\line(0,1){10}}\put(0,10){\circle*{2.5}}\end{picture}}

\newcommand{\inc}{\hookrightarrow}
\newcommand{\onto}{\twoheadrightarrow} 

\newcommand{\QED}{\qed}

\title[Hopf algebra of trees]{Structure of the Loday-Ronco\\ Hopf algebra of trees} 
\author{Marcelo Aguiar}
\address{Department of Mathematics\\
         Texas A\&M University\\
         College Station\\
         TX \ 77843\\
         USA}
\email{maguiar@math.tamu.edu}
\urladdr{http://www.math.tamu.edu/$\sim$maguiar}

\author{Frank Sottile}
\address{Department of Mathematics\\
         Texas A\&M University\\
         College Station\\
         TX \ 77843\\
         USA}
\email{sottile@math.tamu.edu}
\urladdr{http://www.math.tamu.edu/$\sim$sottile}

\thanks{Aguiar supported in part by NSF grant DMS-0302423. 
Sottile supported in part by NSF CAREER
  grant DMS-0134860 and the Clay Mathematical Institute}
\begin{document}

\begin{abstract}
 Loday and Ronco defined an interesting  Hopf algebra structure on the linear
 span of the set of planar binary trees.
 They  showed that the inclusion of the Hopf algebra of non-commutative
 symmetric functions in the Malvenuto-Reutenauer Hopf
 algebra of permutations factors through their Hopf algebra of trees, and these
 maps correspond to natural maps from the weak order on the symmetric group to
 the Tamari order on planar binary trees to the boolean algebra. 

 We further study the structure of this Hopf algebra of trees using a
 new basis for it.
 We describe the product, coproduct, and antipode in terms of this basis
 and use these results to elucidate its Hopf-algebraic structure.
 We also obtain a transparent proof of its isomorphism with the
 non-commutative Connes-Kreimer Hopf algebra of Foissy,
 and show that this algebra is related to non-commutative symmetric functions as
 the (commutative) Connes-Kreimer Hopf algebra is related to symmetric
 functions.  
\end{abstract}

\maketitle
%%%%%%%%%%%%%%%%%%%%%%%%%%%%%%%%%%%%%%%%%%%%%%%%%%%%%%%%%%%%%%%%%%%%%%%%%%%%%%%%%%
%\tableofcontents

\contentsline {section}{\tocsection {}{}{Introduction}}{1}
\contentsline {section}{\tocsection {}{1}{Basic Definitions}}{2}
\contentsline {section}{\tocsection {}{2}{Some Galois connections between posets}}{10}
\contentsline {section}{\tocsection {}{3}{Some Hopf morphisms involving $\YSym$}}{13}
\contentsline {section}{\tocsection {}{4}{Geometric interpretation of the
                                          product of $\YSym$}}{14} 
\contentsline {section}{\tocsection {}{5}{Cofreeness and the coalgebra structure
                                          of $\YSym$}}{17} 
\contentsline {section}{\tocsection {}{6}{Antipode of $\YSym$}}{19}
\contentsline {section}{\tocsection {}{7}{Crossed product decompositions for $\SSym$ and $\YSym$}}{22}
\contentsline {section}{\tocsection {}{8}{The dual of $\YSym$ and the non-commutative Connes-Kreimer Hopf algebra}}{23}
\contentsline {section}{\tocsection {}{}{References}}{31}

%% May uncomment the above for manual table of contents (without subsections)
%%%%%%%%%%%%%%%%%%%%%%%%%%%%%%%%%%%%%%%%%%%%%%%%%%%%%%%%%%%%%%%%%%%%%%%%%%%%%%%%%%
        
\section*{Introduction}
In 1998, Loday and Ronco defined a Hopf algebra $\LR$ on the linear
span of the set of rooted planar binary trees~\cite{LR98}.
This Hopf algebra is the
free dendriform algebra on one generator~\cite{Lod01}.
In~\cite{LR98,LR02}, Loday and Ronco showed how natural poset maps between the
weak order on the symmetric groups, the Tamari order on rooted planar binary
trees with $n$ leaves, and the Boolean posets induce injections of Hopf
algebras
\[
  \NSym\ \hookrightarrow\ \LR\ \hookrightarrow\ \SSym\,,
\]
where $\NSym$ is the Hopf algebra of non-commutative symmetric functions~\cite{GKal}
and $\SSym$ is the Malvenuto-Reutenauer Hopf algebra of
permutations~\cite{MR95}.

Simultaneously, Hopf algebras of trees were proposed by Connes and
Kreimer~\cite{CK98,K98} and Brouder and Frabetti~\cite{BF01,BF03} to encode
renormalization in quantum field theories.
The obvious importance of these algebras led to intense study, and by work of 
Foissy~\cite{Foi02a,Foi02b}, Hivert-Novelli-Thibon~\cite{HNT03},
Holtkamp~\cite{Ho03}, and Van der Laan~\cite{Laan}, the Hopf algebras of
Loday-Ronco, Brouder-Frabetti, and the non-commutative Connes-Kreimer  Hopf
algebra are known to be isomorphic, self-dual, free (and cofree).  

We described the elementary structure of $\SSym$ with respect
to a new basis and used those results to further elucidate its structure as a
Hopf algebra~\cite{AS02}.
Here, we use a similar approach to study $\YSym:=(\LR)^*$, the graded dual Hopf algebra to $\LR$.
We define a new basis for $\YSym$ related to the (dual of) the
Loday-Ronco basis via M\"obius inversion on the poset of trees.
We next describe the elementary structure of $\YSym$ with respect to this new
basis, use those results to show that it is cofree, and then study its relation
to $\SSym$ and $\QSym$, the Hopf algebra of quasi-symmetric functions.
This basis allows us to give an explicit isomorphism between the Loday-Ronco Hopf algebra $\LR$ and the noncommutative Connes-Kreimer Hopf
algebra of Foissy (this coincides with the isomorphism constructed by
Holtkamp~\cite{Ho03} and Palacios~\cite{Pa02}).  We use it to show that a canonical involution of $\QSym$ can be
lifted to $\YSym$ and
to deduce a commutative diagram involving
the Connes-Kreimer Hopf algebras (commutative and non-commutative) on one hand, and symmetric and non-commutative symmetric functions on the other.

Our approach provides a unified framework to understand the
structures of $\YSym$ and explain them in the context of the well-understood
Hopf algebras $\SSym$, $\QSym$, and $\NSym$ of algebraic combinatorics.
A similarly unified approach, through realizations of the algebras 
$\SSym$ and $\YSym$  via
combinatorial monoids, has been recently obtained by Hivert, Novelli, and
Thibon~\cite{HNT02,HNT03,HNT04}. Another interesting approach, involving lattice
congruences, has been proposed by Reading~\cite{Re04a,Re04b}.

%%%%%%%%%%%%%%%%%%%%%%%%%%%%%%%%%%%%%%%%%%%%%%%%%%%%%%%%%%%%%%%%%%%
\section{Basic definitions}\label{S:basic}

\subsection{Compositions, permutations, and trees}\label{S:com-per-tre}
Throughout, $n$ is a non-negative integer and $[n]$ denotes the set
$\{1,2,\ldots,n\}$.
A {\it composition} $\alpha$ of $n$ is a sequence
$\alpha=(\alpha_1,\ldots,\alpha_k)$ of positive integers whose sum is $n$.
Associating the set $I(\alpha):=
\{\alpha_1,\alpha_1+\alpha_2,\ldots,\alpha_1+\cdots+\alpha_{k-1}\}$
to a composition $\alpha$ of $n$ gives a bijection between compositions of $n$
and subsets of $[n{-}1]$. 
Compositions of $n$ are partially ordered by \emph{refinement},
which is defined by its cover relations
 \[
   (\alpha_1,\ldots,\alpha_i+\alpha_{i+1},\ldots,\alpha_k)\ \lessdot\ 
   (\alpha_1,\ldots,\alpha_k)\,.
 \]
Under the association $\alpha\leftrightarrow I(\alpha)$, refinement
corresponds to set inclusion, so we simply identify the poset of compositions
of $n$ with the Boolean poset $\calQ_n$ of subsets of $[n{-}1]$.\smallskip

Let $\frakS_n$ be the group of permutations of $[n]$.
We use one-line notation for permutations, writing 
$\sigma=(\sigma(1),\sigma(2),\ldots,\sigma(n))$
 and sometimes omitting the parentheses and commas.
The {\it standard permutation} $\st(a_1,\ldots,a_p)\in\frakS_p$ of a sequence   
$(a_1,\ldots,a_p)$ of distinct integers is the unique permutation $\sigma$ such that   
\begin{equation*}%\label{E:st}
   \sigma(i)<\sigma(j)\  \iff\  a_i<a_j.
\end{equation*}

An \emph{inversion} in a permutation $\sigma\in\frakS_n$ is a pair of positions
$1\leq i<j\leq n$ with $\sigma(i)>\sigma(j)$.
Let $\Inv(\sigma)$ denote the set of inversions of $\sigma$.
Given $\sigma,\tau\in\frakS_n$, we write $\sigma\leq\tau$ if 
$\Inv(\sigma)\subseteq\Inv(\tau)$.
This defines the (left) \emph{weak order} on $\frakS_n$.
%The length of a permutation $\sigma$ is $\ell(\sigma)=\#\Inv(\sigma)$.
The identity permutation $\id_n$ is the minimum element in 
$\frakS_n$ and $\omega_n=(n,\dotsc,2,1)$ is the maximum. See~\cite[Figure 1]{AS02} for a picture of the weak order on $\frakS_4$.

The grafting of two permutations $\sigma\in\frakS_p$ and $\tau\in\frakS_q$ is the permutation $\sigma\vee\tau\in\frakS_{p+q+1}$  with values
\begin{equation}\label{E:def-grafting-perm}
  \sigma(1){+}q,\ \sigma(2)+q,\ \dotsc,\ \sigma(p){+}q,
   \ \ p{+}q{+}1,\ \ 
   \tau(1),\ \tau(2),\ \dotsc,\ \tau(q)\,.
\end{equation}
Similarly, let $\sigma\underline{\vee}\tau\in\frakS_{p+q+1}$ be the permutation with values
\begin{equation*}
  \sigma(1),\ \sigma(2),\ \dotsc,\ \sigma(p),
   \ \ p{+}q{+}1,\ \ 
   \tau(1){+}p,\ \tau(2){+}p,\ \dotsc,\ \tau(q){+}p\,.
\end{equation*}
This is the operation considered in~\cite[Def. 1.6]{LR02}. As in~\cite[Def. 1.9]{LR02}, let
$\sigma\backslash\tau\in\frakS_{p+q}$ be the permutation whose values are
 \[
   \sigma(1)+q,\, \sigma(2)+q,\,\dotsc,\,\sigma(p)+q,\ \,
     \tau(1),\,\tau(2),\,\ldots,\,\tau(q)\ .
 \]
Consider decompositions $\rho=\sigma\backslash\tau$ of a permutation $\rho$.
Every permutation $\rho$ may be written as 
$\rho=\rho\backslash\id_0 = \id_0\backslash\rho$, where 
$\id_0\in\frakS_0$ is the empty permutation.
A permutation $\rho\neq\id_0$ has \emph{no global descents} if these are its only such
decompositions.
The operation $\backslash$ is associative and a permutation $\rho\neq\id_0$ has a unique decomposition into  permutations with no global descents.
We similarly have the associative operation $/$ to form the permutation
$\sigma/\tau$ (often denoted $\sigma\times\tau$ in the literature) whose values are 
 \[
   \sigma(1),\, \sigma(2),\,\dotsc,\,\sigma(p),\ \;
     \tau(1)+p,\,\tau(2)+p,\,\ldots,\,\tau(q)+p\,.
 \]
The following properties are immediate from the definitions. For any permutations $\rho,\sigma,\tau$,
\begin{align}\label{E:assoc-perm-1}
(\rho\vee\sigma)\backslash\tau &=\rho\vee(\sigma\backslash\tau)\,,\\
\label{E:assoc-perm-2}
\rho/(\sigma\underline{\vee}\tau) &=(\rho/\sigma)\underline{\vee}\tau\,.
\end{align}

\begin{figure}[htb]
\[
  \begin{picture}(320,227) \label{F:Trees}
   \put(0,44){\epsffile{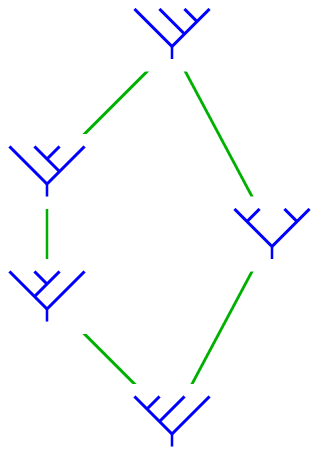}}
   %fig2dev -Leps -m0.8 Y4.fig Y4.eps
   \put(120,0){  
  \begin{picture}(195,220)(-6,-5)
   \put(0,2){\epsfysize=200pt\epsffile{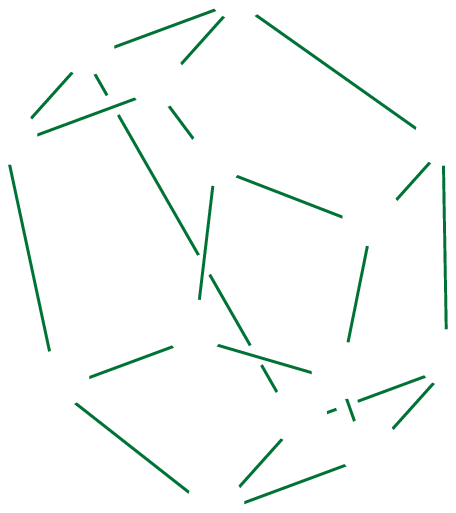}}
   \put(134, 18.4){\epsfysize=14pt\epsffile{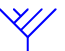}}
   \put(165, 52.6){\epsfysize=14pt\epsffile{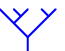}}
   \put(123.2, 48.2){\epsfysize=14pt\epsffile{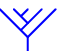}}
   \put(135.4,108.3){\epsfysize=14pt\epsffile{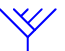}}
   \put(163.4,139.2){\epsfysize=14pt\epsffile{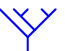}}
   \put( 74.1,131.5){\epsfysize=14pt\epsffile{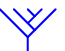}}
   \put( 12.7, 45){\epsfysize=14pt\epsffile{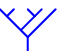}}
   \put( 74.4, -3.6){\epsfysize=14pt\epsffile{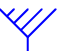}}
   \put( 66.2, 64.5){\epsfysize=14pt\epsffile{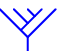}}
   \put( 23.5,174.8){\epsfysize=14pt\epsffile{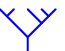}}
   \put( 51.5,162.3){\epsfysize=14pt\epsffile{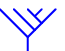}}
   \put( -7.6,140){\epsfysize=14pt\epsffile{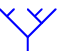}}
   \put( 82.7,196.5){\epsfysize=14pt\epsffile{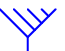}}
   \put(105.6,33.7){\epsfysize=14pt\epsffile{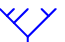}}
  \end{picture}}
  \end{picture}
\]
\caption{The Tamari order on $\calY_3$ and $\calY_4$.} 
\end{figure}

Let $\calY_n$ be the set of rooted, planar binary trees with $n$ interior
nodes (and thus $n+1$ leaves).
The {\em Tamari order}  on $\calY_n$ is the partial order whose cover relations are obtained by moving a child node
directly above a given node from the left to the right branch above the
given node. 
Thus
 \[
  \raisebox{-7pt}{\epsffile{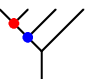}}
      {\color{red} \ \longrightarrow\  }
  \raisebox{-7pt}{\epsffile{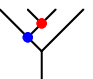}}
    {\color{blue}   \ \longrightarrow\  }
  \raisebox{-7pt}{\epsffile{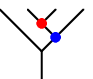}}
      {\color{red} \ \longrightarrow\   }
  \raisebox{-7pt}{\epsffile{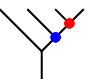}}
 \]
is an increasing chain in $\calY_3$
(the moving vertices are marked with  dots).
Only basic properties of the Tamari order are needed in this paper; their proofs will be provided. For more properties, see~\cite[Sec. 9]{BW97}.
  Figure~\ref{F:Trees} shows the Tamari order on $\calY_3$ and $\calY_4$. 

%%%%%%%%%%%%%%%%%%%%%%%%%%%%%%%%%%%%%%%%%%%%%%%%%%%%%%%%%%%%%%%%%%%%%%%%%
%%%%%%%%%%%%%%%%%%%%%%%%%%%%%%%%%%%%%%%%%%%%%%%%%%%%%%%%%%%

Let $1_n$ be the minimum  tree in $\calY_n$. It is called a {\it right comb} as all of its leaves are right pointing:
\[
  1_4\ =\ \epsffile{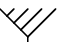}\qquad
  1_7\ =\ \raisebox{-5pt}{\epsffile{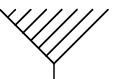}}.
\]
Given trees $s\in\calY_p$ and $t\in\calY_q$, the tree
$s\vee t\in\calY_{p+q+1}$ is obtained by grafting the root of $s$ onto the
left leaf of the tree\ \epsffile{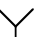}\ and
the root of $t$ onto its right leaf. Below we display trees $s$, $t$, and $s\vee t$,
indicating the position of the grafts with dots.
 \begin{equation*} 
   \epsffile{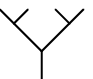}\qquad
 \epsffile{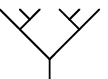}\qquad    
    \epsffile{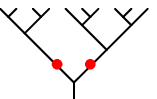}
 \end{equation*}

For $n>0$, every tree $t\in\calY_n$ has a unique decomposition $t=t_l\vee t_r$ with
$t_l\in\calY_p$, $t_r\in\calY_q$, and $n=p+q+1$.
Thus $\calY_n$ is in bijection with
$\bigsqcup_{p+q=n{-}1}\calY_p\times\calY_q$, and
since $\calY_0=\{\,\raisebox{-2pt}{\rule{0.4pt}{11pt}}\,\}$ and
$\calY_1=\{\,\epsffile{figures/1.eps}\}$, we see that 
$\calY_n$ contains the Catalan number $\frac{(2n)!}{n!(n+1)!}$ of
trees. 

For trees $s$ and $t$, let $s\backslash t$ be the tree obtained by adjoining
the root of $t$ to the rightmost branch of $s$.
Similarly, $s/t$ is obtained by grafting the root of $s$ to the leftmost
branch of $t$.
These operations are associative.
Here, we display trees $s$, $t$, $s\backslash t$, and $s/t$,
indicating the position of the graft with a dot.
 \begin{equation*} %\label{E:irred_decomp}
    \epsffile{figures/231-big.eps}\qquad
    \epsffile{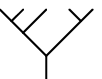}\qquad
    \epsffile{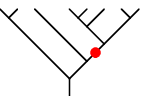}\qquad
    \epsffile{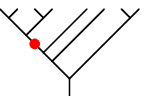}
 \end{equation*}
The following properties are immediate from the definitions. For any trees $s$
and $t$,
 \begin{eqnarray}\label{E:assoc-tree-1}
  s\backslash t & =&s_l\vee(s_r\backslash t)\,,\\
  \label{E:assoc-tree-2}
  s/ t &  =&(s/ t_l)\vee t_r\,.
\end{eqnarray}

A (right) decomposition of a tree $t$ is a way of writing $t$ as 
$r\backslash s$.
Note that $t=t\backslash 1_0= 1_0\backslash t$,
 so every tree has two trivial decompositions.
We say that a tree $t\neq 1_0$ is \emph{progressive} if these are its only right
decompositions. 
For any tree $t\in\calY_n$ we have, by~\eqref{E:assoc-tree-1}, $t=t_l\vee
t_r=t_l\vee(1_0\backslash t_r)=(t_l\vee 1_0)\backslash t_r$. 
Also, for any trees $s,r$ we have $s\backslash r=(s_l\vee s_r)\backslash
r=s_l\vee(s_r \backslash r)$. 
Therefore, $t$ is progressive if and only if $t_r=1_0=|$.
Geometrically, progressive trees have no branching along the right branch from
the root; equivalently, all internal nodes are to the {\em left}  of the root.

Every tree $t\neq 1_0$ has a unique decomposition into progressive trees, 
$t=t_1\backslash t_2 \backslash \cdots \backslash t_k$.
For example,
 \[
  \begin{picture}(220,25)
   \put(0,0){\epsffile{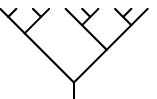}}
   \put(50,10){=}
   \put(70,2.5){\epsffile{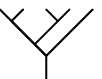}\epsffile{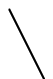}
   \epsffile{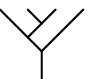}\epsffile{figures/BSL.eps}
   \epsffile{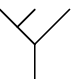}}
   \put(190,10){.}
   \end{picture}
 \]

%%%%%%%%%%%%%%%%%%%%%%%%%%%%%%%%%%%%%%%%%%%%%%%%%%%%%%%%%%%%%%%%%%%%%%%%%%%
\subsection{Some maps of posets}\label{S:combinatorics} 

Order-preserving maps (poset maps) between the posets $\calQ_n$, $\calY_n$,
and $\frakS_n$ are central to the structures of the Hopf algebras 
$\QSym$, $\YSym$, and $\SSym$.
A permutation $\sigma\in\frakS_n$ has a {\it descent} at a position
$p$ if $(p,p{+}1)\in\Inv(\sigma)$, that is if $\sigma(p)>\sigma(p{+}1)$.
Let $\Des(\sigma)\in\calQ_n$ denote the set of descents of a permutation $\sigma$.
Then $\Des\colon\frakS_n\to\calQ_n$ is a surjection of posets.
Given $\setS=\{p_1,\dotsc,p_k\}\in\calQ_n$, let $Z(\setS)\in\frakS_n$ be
 \[
   Z(\setS)\ :=\ \id_{p_1}\,\backslash\, \id_{p_2-p_1}\,\backslash\,
                      \dotsb\,\backslash\, \id_{n-p_k}\,.
 \]
This is the maximum permutation in $\frakS_n$ whose descent set is $\setS$.
The map $Z\colon\calQ_n\hookrightarrow\frakS_n$ is an embedding of posets,
in the sense that $\setS\subseteq\setT\iff Z(\setS)\leq Z(\setT)$. 

The image of $Z$ is described as follows.
A permutation $\sigma\in\frakS_n$ is {\em $132$-avoiding} if whenever $i<j<k\leq n$,
then we {\it do not} have $\sigma(i)<\sigma(k)<\sigma(j)$.
For example, $43512$ is $132$-avoiding. Similarly, $\sigma$ is $213$-avoiding if
whenever $i<j<k\leq n$,
then we  do not have $\sigma(j)<\sigma(i)<\sigma(k)$.
The definition of $Z$ implies that $Z(\setS)$ is
both $132$ and $213$-avoiding.
Since the number of $(132,213)$-avoiding permutations is $2^{n-1}$~\cite[Ch. 14, Ex. 4]{Bo03}, the map
$Z$ embeds 
$\calQ_n$ as the subposet of $\frakS_n$ consisting of $(132,213)$-avoiding
permutations.  

There is a well-known map  that sends a 
permutation  to a tree~\cite[pp.~23-24]{St86},~\cite[Def. 9.9]{BW97}. 
We are interested in the
following variant $\lambda\colon\frakS_n\to\calY_n$, as considered
in~\cite[Section 2.4]{LR98}.  
We define $\lambda(\id_0)=1_0$. For $n\geq 1$, let $\sigma\in\frakS_n$
and  $j:=\sigma^{-1}(n)$.
We set $\sigma_l:=\st(\sigma(1),\ldots,\sigma(j{-}1))$,
$\sigma_r:=\st(\sigma(j{+}1),\ldots,\sigma(n))$, and define 
 \begin{equation}\label{E:def-lambda}
  \lambda(\sigma)\ :=\ \lambda(\sigma_l)\vee\lambda(\sigma_r)\,.
 \end{equation} 
In other words, we construct $\lambda(\sigma)$ recursively by grafting
$\lambda(\sigma_l)$ and $\lambda(\sigma_r)$ onto the left and right branches of
$\,$\epsffile{figures/1.eps}. 
For example, if $\sigma=\Red{564}9\Blue{73812}$ then $j=4$, $\sigma_l=\Red{231}$, 
$\sigma_r=\Blue{43512}$,  and
\[
   \begin{picture}(360,45)(-84,0)
    \put(-84,18){$\lambda(\sigma_l)=$}
    \put(-42,12){\epsffile{figures/231-big.eps}}
     \put(0,18){$\lambda(\sigma_r)=$}
     \put(40,12){\epsffile{figures/43512-big.eps}}
     \put(85,18){$\Longrightarrow$}
  \put(120,18){$\lambda(\sigma)$}
    \put(150,18){=}
    \put(158,20){\epsffile{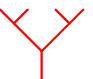}}
      \put(170,-3){\epsffile{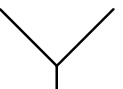}}
   \put(188.4,20){\epsffile{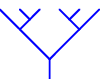}}
     \put(217,18){=}
     \put(230,7){\epsffile{figures/786943512.eps}}       
\end{picture}
 \]
Note that if $\sigma\in\frakS_p$ and $\tau\in\frakS_q$ then
$\lambda(\sigma\vee\tau)=\lambda(\sigma)\vee\lambda(\tau)$. 

It is known that $\lambda$ is a surjective morphism of
posets~\cite[Prop. 9.10]{BW97},~\cite[Cor. 2.8]{LR02}. 

Consider the maps $\gamma,\underline{\gamma}\colon\calY_n\to\frakS_n$ 
defined recursively by $\gamma(1_0)=\underline{\gamma}(1_0):=\id_0$ and
 \begin{equation}\label{E:def-gamma}
  \gamma(t)\ :=\ \gamma(t_l)\vee\gamma(t_r) \qquad\text{and}\qquad
  \underline{\gamma}(t)\ :=\ \underline{\gamma}(t_l)\underline{\vee}\underline{\gamma}(t_r)\,.
 \end{equation} 
 These are the maps denoted $\mathrm{Max}$ and $\mathrm{Min}$ by Loday and Ronco~\cite[Def. 2.4]{LR02}. They  show that~\cite[Thm. 2.5]{LR02}
 \begin{equation}\label{E:gamma-lambda}
  \gamma(t)\ :=\ \max\{ \sigma\in\frakS_n\mid \lambda(\sigma)=t\}\qquad\text{and}\qquad
  \underline{\gamma}(t)\ :=\ \min\{ \sigma\in\frakS_n\mid \lambda(\sigma)=t\}\,.
 \end{equation}
In particular, both $\gamma$ and $\underline{\gamma}$ are sections of $\lambda$. In this paper, we are mostly concerned with the map $\gamma$.
The recursive definition of $\gamma$ implies that $\gamma(t)$ is
132-avoiding.
Since $\calY_n$ and the set of $132$-avoiding permutations in $\frakS_n$ are 
equinumerous~\cite[p.~261]{St99}, the map $\gamma$ embeds
$\calY_n$ as the subposet of $\frakS_n$ consisting of $132$-avoiding
permutations. 

Since $Z(\setS)$ is 132-avoiding,  there is a unique
map $C:\calQ_n\hookrightarrow\calY_n$ such that $Z=\gamma\circ C$.
It follows that $C$ is an embedding of posets. Explicitly, if $\setS=\{p_1,\dotsc,p_k\}\in\calQ_n$, then
\begin{equation}\label{E:def-C}
   C(\setS)\ :=\ 1_{p_1}\,\backslash\, 1_{p_2-p_1}\,\backslash\,
                      \dotsb\,\backslash\, 1_{n-p_k}\,.
\end{equation}

A tree $t\in\calY_n$ has $n{+}1$ leaves, which we number from $1$ to $n{-}1$
left-to-right, excluding the two outermost leaves. 
Let $L(t)$ be the set of labels of those leaves that  point left. 
For the tree $t\in\calY_8$ below, $L(t)=\{2,5,7\}$.
\[
  \begin{picture}(100,65)
   \put(0,0){\scalebox{2}{\epsffile{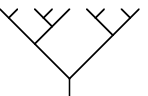}}}
   \put(7,55){1}\put(27,55){3}\put(37,55){4}
   \put(57,55){6}
   {\color{red}
   \put(17,55){2}\put(47,55){5}\put(67,55){7}}
 \end{picture}
\]
Loday and Ronco~\cite[Sec. 4.4]{LR98} note that $\Des=\lambda\circ L$. It
follows that 
$L(t)=\Des(\gamma(t))$ and $L$ is a surjective morphism of posets. In summary:

\begin{prop}
 We have the following commutative diagrams of poset maps.
 \begin{equation*}%\label{E:inclusions}
   \begin{picture}(113,36)(0,-1)

          \put(48,28){$\calY_n$}
    \put(0,0){$\calQ_n$}    \put(97,0){$\frakS_n$}
    \put(45,27){\vector(-2,-1){30}} \put(45,27){\vector(-2,-1){25}}
    \put(92.5, 2){\vector(-1,0){75}}\put(92.5, 2){\vector(-1,0){70}}
    \put(95,12){\vector(-2,1){30}}\put(95,12){\vector(-2,1){25}}
  
    \put(23,23){\scriptsize$L$}    \put(82,22){\scriptsize$\lambda$}
          \put(48,5){\scriptsize$\Des$}

   \end{picture} 
  \qquad
   \begin{picture}(113,36)(0,-1)

          \put(48,28){$\calY_n$}
    \put(0,0){$\calQ_n$}    \put(97,0){$\frakS_n$}
    \put(15,12){\epsffile{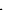}}
    \put(15,12){\vector(2,1){30}}
    \put(18, 2){\epsffile{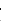}}
    \put(18, 2){\vector(1,0){75}}
    \put(65,25){\epsffile{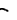}}
    \put(65,25){\vector(2,-1){30}}
  
    \put(24,24){\scriptsize$C$}    \put(78,22){\scriptsize$\gamma$}
          \put(51,6){\scriptsize$Z$}

   \end{picture} 
 \end{equation*}
 In addition, $\Des\circ Z=\id_{\calQ_n}$, $\lambda\circ\gamma=\id_{\calY_n}$,
 $L\circ C=\id_{\calQ_n}$.
\end{prop}

The maps $\lambda$, $\gamma$, and $\underline{\gamma}$ are well-behaved with respect to the operations
$\backslash$ and $/$.

\begin{prop}\label{P:starslash}
 Let $\sigma\in\frakS_p$,  $\tau\in\frakS_q$, $s\in\calY_p$, and $t\in\calY_q$.
 Then
 \begin{eqnarray}
 \label{E:lambda-slash}
   \lambda(\sigma\backslash\tau)\ =\ \lambda(\sigma)\backslash\lambda(\tau)
     &\quad\textrm{and}\quad&  
    \lambda(\sigma/\tau)\ =\ \lambda(\sigma)/\lambda(\tau) \\
 \label{E:gamma-slash}
 \gamma(s\backslash t)\ =\ \gamma(s)\backslash\gamma(t) &\quad\textrm{and}\quad&  \underline{\gamma}(s/t)\ =\ \underline{\gamma}(s)/\underline{\gamma}(t)\,.
 \end{eqnarray}
\end{prop}
\noindent{\it Proof. } 
The assertions about $\lambda$ are given in \cite[Thm. 2.9]{LR02}.
For $\gamma$, note that by \eqref{E:assoc-tree-1}, 
$s\backslash t  =s_l\vee(s_r\backslash t)$. 
Since by definition~\eqref{E:def-gamma} $\gamma$ preserves the grafting
operations $\vee$, we have $\gamma(s\backslash
t)=\gamma(s_l)\vee\gamma(s_r\backslash t)$. 
Proceeding inductively we derive $\gamma(s\backslash
t)=\gamma(s_l)\vee\bigl(\gamma(s_r)\backslash \gamma(t)\bigr)$. 
Finally, from~\eqref{E:assoc-perm-1} we conclude 
$\gamma(s\backslash t)=\bigl(\gamma(s_l)\vee\gamma(s_r)\bigr)\backslash \gamma(t)=\gamma(s)\backslash\gamma(t)$. The assertion about $\underline{\gamma}$ can be similarly obtained from~\eqref{E:assoc-tree-2} and~\eqref{E:assoc-perm-2}.
\quad\QED\medskip

We discuss these maps further in Section~\ref{S:galois}.

\begin{rem}
 The weak order on $\frakS_n$ was defined by the inclusion of inversion
 sets.
 If we define the inversion set of a tree $t$ to be the inversion set of the permutation $\gamma(t)$, then we obtain
 that for trees $s,t\in\calY_n$, $s\leq t \Leftrightarrow \Inv(s)\subseteq \Inv(t)$.
 Hugh Thomas pointed out that this inversion set can de described directly in terms of the tree as follows. 
 Suppose that planar binary trees are drawn with branches either pointing
 right (\raisebox{-1pt}{\includegraphics{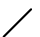}}) or left
 (\raisebox{-1pt}{\includegraphics{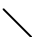}}). 
 If we illuminate a tree $t$ from the Northeast, then its inversion set
 $\Inv(t)$ is the shaded region among its branches.
 For example, Figure~\ref{fig:treeInversion} shows the inversion set for the
 tree $\lambda(41253)$.
 \begin{figure}[htb]
   \epsfysize=100pt\epsffile{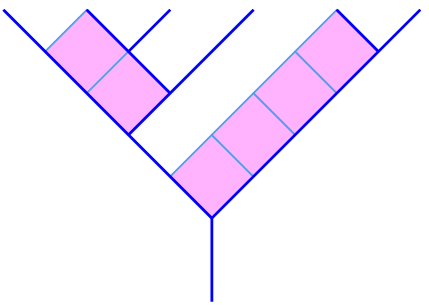} \qquad
   \epsfysize=100pt\epsffile{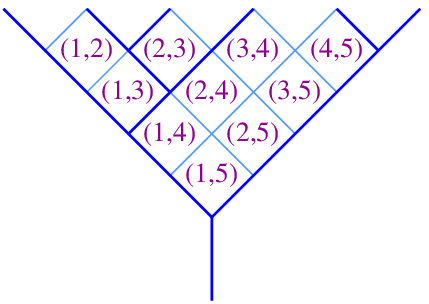}
 \caption{The inversion set of a tree. \label{fig:treeInversion}}
 \end{figure}
 We relate this to the inversion set of the permutation $\gamma(t)$.
 If we draw trees in $\calY_n$ on a tilted grid rotated $45^\circ$ counterclockwise (as we do), then the region between the leftmost and
 rightmost branches is a collection of boxes which may be labeled from $(1,2)$
 in the leftmost box to $(n{-}1,n)$ in the rightmost box, with the first index
 increasing in the Northeast direction and the second in the Southeast
 direction.   (See  Figure~\ref{fig:treeInversion}.)
 Then, given a tree $t$,  the labels of the boxes in the shade is the inversion
 set  of $\gamma(t)$.
 For example, if $t=\lambda(41253)$, then $\gamma(t)=42351$, and we see that
\[
   \Inv(t)\ =\ \{(1,2),\,(1,3),\,(1,5),\,(2,5),\,(3,5),\,(4,5)\}\ =\ 
   \Inv(42351)\,.
\]
\end{rem}

%%%%%%%%%%%%%%%%%%%%%%%%%%%%%%%%%%%%%%%%%%%%%%%%%%%%%%%%%%%%%%%%%%%%%%%%%%%%%%%
\subsection{The Hopf algebra of quasi-symmetric functions}\label{S:hopfquasi}
We use elementary properties of Hopf algebras, as given in 
the book~\cite{Mo93a}.
Our Hopf algebras $H$ will be graded connected Hopf algebras over $\Q$.
Thus the $\Q$-algebra $H$ is the direct sum 
$\bigoplus\{H_n\mid n=0,1,\ldots\}$ of its 
homogeneous components $H_n$, with $H_0=\Q$, the product and
coproduct respect the grading, and the counit is the projection onto $H_0$.

The algebra $\QSym$ of quasi-sym\-met\-ric functions was introduced by
 Gessel~\cite{Ges} in connection to work of Stanley~\cite{St72}. Malvenuto
described its Hopf algebra structure~\cite[Section 4.1]{Malv}. 
See also~\cite[9.4]{Re93} or~\cite[Section 7.19]{St99}.

$\QSym$ is a graded connected Hopf algebra. The component of degree $n$ has a linear basis of {\em monomial
quasi-symmetric functions} $\setM_\alpha$ indexed
by compositions $\alpha=(a_1,\ldots,a_k)$ of $n$. The coproduct is
\begin{equation}\label{E:coproduct-QSym}
\Delta(M_\alpha)=\sum_{i=0}^k M_{(a_1,\ldots,a_i)}\otimes M_{(a_{i+1},\ldots,a_k)}\,.
\end{equation}
The product of two monomial functions $M_\alpha$ and $M_\beta$ can be described in terms of quasi-shuffles of $\alpha$ and $\beta$. A geometric
description for the structure constants in terms of faces of the cube was given in~\cite[Thm. 7.6]{AS02}.

Gessel's {\em fundamental}  quasi-symmetric function $\setF_\beta$ is defined by
 \[
    \setF_\alpha\ =\ \sum_{\alpha\leq\beta} \setM_\beta\,.
 \]
By M\"obius inversion, we have
 \[
    \setM_\alpha\ =\ \sum_{\alpha\leq\beta}(-1)^{k(\beta)-k(\alpha)}\setF_\beta\,,
 \]
where $k(\alpha)$ is the number of parts of $\alpha$. Thus the set $\{\setF_\alpha\}$ forms another basis of $\QSym$. 

We often index these monomial and fundamental quasi-symmetric functions by
subsets of $[n{-}1]$.
Accordingly, given a composition $\alpha$ of $n$ with $\setS=I(\alpha)$, we
define
 \[
    \setF_{\setS}\ :=\ \setF_\alpha \qquad\text{and}\qquad
    \setM_{\setS}\ :=\ \setM_\alpha\,.
 \]
This notation suppresses the dependence on $n$, which is 
understood from the context.

%%%%%%%%%%%%%%%%%%%%%%%%%%%%%%%%%%%%%%%%%%%%%%%%%%%%%%%%%%%%%%%%%%%%
\subsection{The  Hopf algebra of permutations}\label{S:MR}

Set $\frakS_\infty:=\bigsqcup_{n\geq 0}\frakS_n$.
Let $\SSym$ be the graded vector space over $\Q$ with {\it fundamental basis}
$\{\calF_\sigma\mid \sigma\in\frakS_\infty\}$, whose degree $n$ component is
spanned by $\{\calF_\sigma\mid \sigma\in\frakS_n\}$. 
Write $1$ for the basis element of degree $0$.
Malvenuto and Reutenauer~\cite{Malv,MR95} described a Hopf algebra structure
on this space that was further elucidated in~\cite{AS02}.
Here, as in~\cite{AS02}, we study the self-dual Hopf algebra $\SSym$ with
respect to bases dual to those in~\cite{Malv,MR95}. 

The product of two basis elements is obtained by shuffling the corresponding
permutations.
For $p,q>0$, set
 \[
   \Sh{(p,q)} \ :=\
              \{\zeta\in \frakS_{p+q}\mid \zeta \mbox{ has at most one
              descent, at position $p$}\}\,.
 \]
This is the  collection of minimal  representatives of left cosets
of $\frakS_p\times\frakS_q$ in $\frakS_{p+q}$. 
These are sometimes called $(p,q)$-shuffles.
For $\sigma\in \frakS_p$ and $\tau\in \frakS_q$, set
 \begin{equation*}%\label{E:prod-fundamental}
  \calF_\sigma \cdot \calF_\tau\ =\ \sum_{\zeta\in \Sh{(p,q)}}
                            \calF_{(\sigma/\tau)\cdot\zeta^{-1}}\,. 
\end{equation*}
This endows $\SSym$ with the structure of a graded algebra with unit 1.

The algebra $\SSym$ is also a graded coalgebra with coproduct given by all
ways of splitting a permutation.
More precisely, define $\Delta\colon\SSym\to\SSym\,\ten\SSym$ by
\begin{equation*} %\label{E:cop-malvenuto}
  \Delta(\calF_\sigma)\ =\ 
   \sum_{p=0}^n \calF_{\st(\sigma(1),\,\ldots,\,\sigma(p))}\ten
                \calF_{\st(\sigma(p{+}1),\,\ldots,\,\sigma(n))}\,,
\end{equation*}
when $\sigma\in\frakS_n$.
With these definitions, $\SSym$ is a graded connected Hopf algebra.

The descent map induces a morphism of Hopf algebras.
 \begin{equation*} %\label{E:descentmap}
   \begin{array}{rcrcl}
     \calD &:& \SSym&\longrightarrow& \QSym\\
           & &\calF_\sigma&\longmapsto& \setF_{\Des(\sigma)}\rule{0pt}{14pt}
   \end{array}
 \end{equation*}

There is another basis $\{\calM_\sigma\mid \sigma\in\frakS_\infty\}$ for
$\SSym$.
For each $n\geq 0$ and $\sigma\in \frakS_n$, define
 \begin{equation}\label{E:def-monomial}
   \calM_\sigma\ :=\ \sum_{\sigma\leq\tau} 
                    \mu_{\frakS_n}(\sigma,\tau)\cdot \calF_v\,,
 \end{equation}
where $\mu_{\frakS_n}(\cdot,\cdot)$ 
is the M\"obius function of the weak order on $\frakS_n$. 
By M\"obius inversion, 
 \begin{equation*}%\label{E:fun-mon}
   \calF_\sigma\ :=\ \sum_{\sigma\leq\tau} \calM_\tau\,,
 \end{equation*}
so these elements $\calM_\sigma$ indeed form a basis of $\SSym$. 
The algebraic structure of $\SSym$ with respect to this $\calM$-basis was
determined in~\cite{AS02}.

\begin{prop}\label{P:MRstuff}
 Let $w\in\frakS_n$.
 Then 
 \begin{align}
 \label{E:MRstuff1}
  \Delta(\calM_\rho) &= \sum_{\rho=\sigma\backslash\tau}
                            \calM_\sigma\ten\calM_\tau\,,\\
 \label{E:MRstuff2}
  \calD(\calM_\sigma) &= \left\{\begin{array}{rcl}
             \setM_\setS&\ &\textrm{if } \sigma=Z(\setS), \ 
                  \mbox{ for some }\ \setS\in\calQ_n\,,\\
                 0      &  &\textrm{otherwise.}
            \end{array}\right.
 \end{align}
\end{prop}

The multiplicative structure constants are non-negative integers
with the following description.
The 1-skeleton of the permutahedron $\Pi_n$ is the Hasse diagram of the
weak order on $\frakS_n$.
Its facets are canonically isomorphic to products of lower
dimensional permutahedra. 
Say that a facet isomorphic to $\Pi_p\times\Pi_q$ has {\it type $(p,q)$}.
Given $\sigma\in\frakS_p$ and $\tau\in\frakS_q$, such a facet has a distinguished
vertex corresponding to $(\sigma,\tau)$ under the canonical isomorphism.
Then, for $\rho\in\frakS_{p+q}$, the coefficient of
$\calM_\rho$ in $\calM_\sigma\cdot\calM_\tau$ is the number of facets of 
$\Pi_{p+q}$ of type $(p,q)$ with the property that the
distinguished vertex is below $\rho$ (in the weak order) and closer to $\rho$ than
to any other vertex in the facet. 

%%%%%%%%%%%%%%%%%%%%%%%%%%%%%%%%%%%%%%%%%%%%%%%%%%%%%%%%%%%%%%%%%%%%%%%%%%%%%%%
\subsection{The Hopf algebra of planar binary trees}\label{S:LR}
Set $\calY_\infty:=\bigsqcup_{n\geq 0}\calY_n$.
The Hopf algebra $\YSym$ of planar binary trees is the graded
vector space over $\Q$ with {\it fundamental basis}
$\{F_t\mid t\in\calY_\infty\}$ whose degree $n$ component is 
spanned by $\{F_t\mid t\in\calY_n\}$.
We describe its multiplication and comultiplication in terms of a
geometric construction on trees.
For any leaf of a tree $t$, we may divide $t$ into two pieces---the
piece left of the leaf and the piece right of the leaf---by dividing
$t$ along the path from the leaf to the root.
We illustrate this on the tree $\lambda(67458231)$.
\[
  %fig2dev -Leps -m0.25 divide.fig divide.eps
  \epsffile{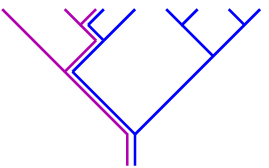}  
    \qquad\raisebox{15pt}{$\Longrightarrow$}\qquad
   %fig2dev -Leps -m0.3 divided.fig divided.eps
   \raisebox{5pt}{\epsffile{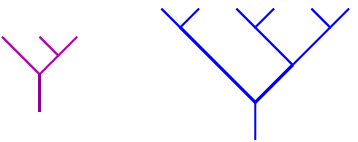}}
\]
If $r$ is the piece to the left of the leaf of $t$ and $s$ the piece
to the right, write $t\to(r,s)$.

Suppose that $t\in\calY_p$ and $s\in\calY_q$. 
Divide  $t$ into $q+1$ pieces at a multisubset of its $p+1$ leaves of
cardinality $q$: 
\[
  t\ \to\ (t_0,t_1,\dotsc,t_q)\,.
\]
This may be done in $\binom{p+q}{p}$ ways. 
Label the leaves of $s$ from $0$ to $q$ left-to-right. 
For each such division $t\to(t_0,t_1,\dotsc,t_q)$, attach $t_i$ to the $i$th leaf of
$s$ to obtain the tree $(t_0,t_1,\dotsc,t_q)/s$.
For example, we divide $\lambda(45231)$ at three leaves to obtain
\[
  %fig2dev -Leps -m0.3 divide45231.fig divide45231.eps
  \epsffile{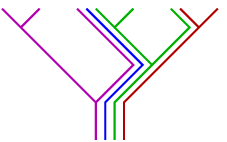}  
    \qquad\raisebox{18pt}{$\rightarrow$}\qquad
   %fig2dev -Leps -m0.3 divided45231.fig divided45231.eps
   \raisebox{17pt}{$\Bigl($}
   \raisebox{10pt}{\begin{picture}(92,23)
    \put(1,0){\,\epsffile{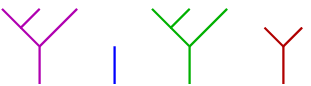}}
    \put(25,1){,} \put(42,1){,} \put(72,1){,}
   \end{picture}}
   \raisebox{17pt}{$\Bigr)$}
    \qquad\raisebox{15pt}{$=\ (t_0,t_1,t_2,t_3)$ .}
\]
Then if $s=\lambda(213)$, the tree $(t_0,t_1,t_2,t_3)/s$ is
\[
  % fig2dev -Leps -m0.3 graft.fig graft.eps
  \epsffile{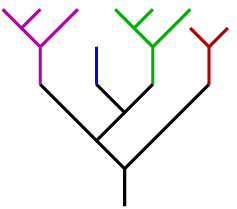}  
   \qquad\raisebox{20pt}{$=$}\qquad
  % fig2dev -Leps -m0.45 grafted.fig grafted.eps
  \raisebox{5pt}{\epsffile{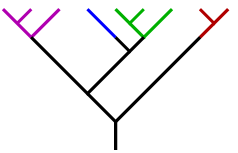}}
\]

We define a coproduct and a product on $\YSym$.
For $t\in\calY_p$ and $s\in\calY_q$, set 
\[
  \Delta (F_t)\ =\ \sum_{t\to(t_0,t_1)} F_{t_0}\otimes F_{t_1}
   \qquad\mbox{and}\qquad
  F_t\cdot F_s\ =\ 
    \sum_{t\to(t_0,t_1,\dotsc,t_q)} F_{(t_0,t_1,\dotsc,t_q)/s}\,.
\]
These are compatible with the operations on $\SSym$ and $\QSym$.
The maps $\lambda\colon \frakS_n\to\calY_n$ and
$L\colon\calY_n\to\calQ_n$ induce linear maps
\begin{equation}\label{E:deflambda}
  \begin{array}{rcrclcrcrcl}
   \Lambda&\colon&\SSym&\longrightarrow&\YSym&\qquad
   \calL  &\colon&\YSym&\longrightarrow&\QSym \\
   &&F_\sigma&\longmapsto& F_{\lambda(\sigma)}&
   &&F_t    &\longmapsto& F_{L(t)}
  \end{array}
\end{equation}

\begin{prop}\label{P:ontoHopf}
 The maps $\Lambda\colon\SSym\to\YSym$ and 
 $\calL\colon\YSym\to\QSym$ are surjective morphisms of Hopf algebras.
\end{prop}

As with both $\QSym$ and $\SSym$, we define another basis 
$\{M_t \mid t\in\calY_\infty\}$, related to the
fundamental basis via M\"obius inversion on $\calY_n$.
For each $n\geq 0$ and $t\in \calY_n$, define
 \begin{equation}\label{E:def-Ymonomial}
   M_t\ :=\ \sum_{t\leq s}\, \mu_{\calY_n}(t,s)\cdot F_s\,,
 \end{equation}
where $\mu_{\calY_n}(\cdot,\cdot)$ 
is the M\"obius function of $\calY_n$.
By M\"obius inversion, 
 \begin{equation*} %\label{E:Yfun-mon}
   F_t\ :=\ \sum_{t\leq s} M_s\,,
 \end{equation*}
so these elements $M_t$ indeed form a basis of $\YSym$. 
For instance,
 \[
   M_{\,\epsffile{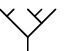}}\ =\ 
   F_{\,\epsffile{figures/3412.eps}} - 
   F_{\,\epsffile{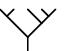}} -
   F_{\,\epsffile{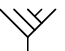}} +
   F_{\,\epsffile{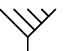}}\,.
 \]

In this paper we determine the algebraic structure of $\YSym$ with respect
to this basis.
For example, we will show (Theorems~\ref{T:maps} and~\ref{T:coproduct}) that,
given $t\in\calY_n$, then  
 \begin{align*}
  \Delta(M_t) & = \sum_{t=r\backslash s} M_r\ten M_s\,;\\
   \Lambda(\calM_\sigma) &= \left\{\begin{array}{lcl}
         M_{t}&\ &\textrm{if $\sigma=\gamma(t)$, $t\in\calY_n$,}\\
            0&&\textrm{otherwise;} \end{array}\right.\\
    \calL( M_t) & = \left\{\begin{array}{ccl}
             \setM_\setT&\ &\textrm{if $t=C(\setT)$,  $\setT\in\calQ_n$,}\\
                 0      &  &\textrm{otherwise.}
            \end{array}\right.
 \end{align*}
We will also obtain a geometric description for the structure constants of the multiplication of $\YSym$ on this basis in terms of the associahedron (Corollary~\ref{C:combprodY}), and an explicit description for the structure constants of the antipode (Theorem~\ref{T:antipodeY}).
%%%%%%%%%%%%%%%%%%%%%%%%%%%%%%%%%%%%%%%%%%%%%%%%%%%%%%%%%%%%%%%%%%%%%%%%%%%%%
\section{Some Galois connections between posets}\label{S:galois}

In Section~\ref{S:combinatorics} we described order-preserving maps
\[
   \begin{picture}(113,36)(0,-1)

          \put(48,28){$\calY_n$}
    \put(0,0){$\calQ_n$}    \put(97,0){$\frakS_n$}
    \put(45,27){\vector(-2,-1){30}} \put(45,27){\vector(-2,-1){25}}
    \put(92.5, 2){\vector(-1,0){75}}\put(92.5, 2){\vector(-1,0){70}}
    \put(95,12){\vector(-2,1){30}}\put(95,12){\vector(-2,1){25}}
  
    \put(23,23){\scriptsize$L$}    \put(82,22){\scriptsize$\lambda$}
          \put(48,5){\scriptsize$\Des$}

   \end{picture} 
  \qquad
   \begin{picture}(113,36)(0,-1)

          \put(48,28){$\calY_n$}
    \put(0,0){$\calQ_n$}    \put(97,0){$\frakS_n$}
    \put(15,12){\epsffile{figures/duhook.eps}}
    \put(15,12){\vector(2,1){30}}
    \put(18, 2){\epsffile{figures/rhook.eps}}
    \put(18, 2){\vector(1,0){75}}
    \put(65,25){\epsffile{figures/ddhook.eps}}
    \put(65,25){\vector(2,-1){30}}
  
    \put(24,24){\scriptsize$C$}    \put(78,22){\scriptsize$\gamma$}
          \put(51,6){\scriptsize$Z$}

   \end{picture} 
\]
between the posets $\calQ_n$, $\calY_n$, and $\frakS_n$. 
Recall from Sections~\ref{S:MR} and~\ref{S:LR} that when the maps in the
leftmost diagram are applied to the fundamental bases,
they induce morphisms of Hopf algebras
$\SSym\twoheadrightarrow\YSym\twoheadrightarrow\QSym$.
The values of these morphisms on the monomial bases can be described through
another set of poset maps given below.
\[
   \begin{picture}(113,36)(0,-1)

          \put(48,28){$\calY_n$}
    \put(0,0){$\calQ_n$}    \put(97,0){$\frakS_n$}
    \put(45,27){\vector(-2,-1){30}} \put(45,27){\vector(-2,-1){25}}
    \put(92.5, 2){\vector(-1,0){75}}\put(92.5, 2){\vector(-1,0){70}}
    \put(95,12){\vector(-2,1){30}}\put(95,12){\vector(-2,1){25}}
  
    \put(23,23){\scriptsize$R$}    \put(82,23){\scriptsize$\rho$}
          \put(51,5){\scriptsize$\GDes$}
   \end{picture} 
\]

A permutation $\sigma$ has a {\it global descent} at a position $p\in[n-1]$ if
$\sigma=\rho\backslash\tau$ with $\rho\in\frakS_p$.
The map $\GDes$ sends a permutation to its set of global descents. 
Global descents were studied in~\cite{AS02}
in connection to the structure of the Hopf algebra $\SSym$.

To define the map $R$, take a tree $t\in\calY_n$ and number its leaves from $1$
to $n{-}1$ left-to-right, excluding the two outermost leaves as before.  
Let $R(t)$ be the set of labels of those leaves that belong to a branch that
emanates from the rightmost branch of the tree. 
In other words, $R(t)$ is set of $j\in [n-1]$ for which the tree 
admits a decomposition $r\backslash t$ with $r\in\calY_j$.
For the tree $t\in\calY_8$ below, $R(t)=\{5,7\}$.
\[ \begin{picture}(100,65)
\put(0,0){\scalebox{2}{\epsffile{figures/67458231.eps}}}
\put(7,55){1}\put(27,55){3}\put(37,55){4}
\put(57,55){6}\put(17,55){2}
{\color{red}
\put(47,55){5}\put(67,55){7}}
\end{picture}\]

The map $\rho$ is defined below. It appears to be new, but by
Theorem~\ref{T:galois} and Lemma~\ref{L:rho} below it is quite natural.

These maps have very interesting order-theoretic properties.
A \emph{Galois connection} between posets $P$ and $Q$ is a pair
$(f,g)$ of order-preserving maps $f\colon P\to Q$ and $g\colon Q\to P$ such
that for any $x\in P$ and $y\in Q$,
 \begin{equation*} %\label{E:galois}
   f(x)\ \leq\ y \ \iff\  x\ \leq\ g(y)\,.
 \end{equation*}
We also say that $f$ is left adjoint to $g$, and $g$ is right adjoint to $f$.

\begin{thm}\label{T:galois}
 We have the following commutative diagrams of order-preserving maps.
\[
   \begin{picture}(113,36)(0,-1)

          \put(48,28){$\calY_n$}
    \put(0,0){$\calQ_n$}    \put(97,0){$\frakS_n$}
    \put(45,27){\vector(-2,-1){30}} \put(45,27){\vector(-2,-1){25}}
    \put(92.5, 2){\vector(-1,0){75}}\put(92.5, 2){\vector(-1,0){70}}
    \put(95,12){\vector(-2,1){30}}\put(95,12){\vector(-2,1){25}}
  
    \put(23,23){\scriptsize$L$}    \put(82,22){\scriptsize$\lambda$}
          \put(48,5){\scriptsize$\Des$}

   \end{picture} 
  \qquad
   \begin{picture}(113,36)(0,-1)

          \put(48,28){$\calY_n$}
    \put(0,0){$\calQ_n$}    \put(97,0){$\frakS_n$}
    \put(15,12){\epsffile{figures/duhook.eps}}
    \put(15,12){\vector(2,1){30}}
    \put(18, 2){\epsffile{figures/rhook.eps}}
    \put(18, 2){\vector(1,0){75}}
    \put(65,25){\epsffile{figures/ddhook.eps}}
    \put(65,25){\vector(2,-1){30}}
  
    \put(24,24){\scriptsize$C$}    \put(78,22){\scriptsize$\gamma$}
          \put(51,6){\scriptsize$Z$}

   \end{picture} 
  \qquad
   \begin{picture}(113,36)(0,-1)

          \put(48,28){$\calY_n$}
    \put(0,0){$\calQ_n$}    \put(97,0){$\frakS_n$}
    \put(45,27){\vector(-2,-1){30}} \put(45,27){\vector(-2,-1){25}}
    \put(92.5, 2){\vector(-1,0){75}}\put(92.5, 2){\vector(-1,0){70}}
    \put(95,12){\vector(-2,1){30}}\put(95,12){\vector(-2,1){25}}
  
    \put(23,23){\scriptsize$R$}    \put(82,23){\scriptsize$\rho$}
          \put(51,5){\scriptsize$\GDes$}
   \end{picture} 
\]
 Moreover, the corresponding maps in adjacent diagrams form Galois connections
 between the appropriate posets.
 That is, the maps in the left diagram are left adjoint to the
 corresponding maps in the central diagram, and the maps in the right diagram
 are right adjoint to the corresponding maps in the central diagram.
\end{thm}

Recall the recursive definition~\eqref{E:def-lambda} of the map 
$\lambda$, where we split a permutation at its greatest value.
The map $\rho$ is similarly described in terms of splitting the permutation,
except now we split it at its first global descent.

We define $\rho(\id_0)=1_0$. For $n\geq 1$, let $\sigma\in\frakS_n$ and suppose
that $j$ is the position of its first global descent.
Let $\sigma_l=\st(\sigma(1),\dotsc,\sigma(j-1))$ and 
$\sigma_r=\st(\sigma(j+1),\dotsc,\sigma(n))$. 
If there are no global descents, we set $j=n$ (and $\sigma_r=\id_0$). 
Note that $\sigma_l$ and $\sigma_r$ generally differ from the permutations
in the definition~\eqref{E:def-lambda} of $\lambda$. 
Define 
\begin{equation}\label{E:def-rho}
\rho(\sigma)\ :=\ \rho(\sigma_l)\vee\rho(\sigma_r)\,.
\end{equation}

For example, if $\sigma=564973812$, then the first global descent occurs at the
position of the 8, and thus
$\sigma_l=\st(564973)=342761$, $\sigma_r=\st(12)=12$, and 
\[
   \begin{picture}(321,45)(-40,-3)
    \put(-84,18){$\rho(\sigma_l)=$}
    \put(-45,12){\epsffile{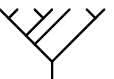}}
     \put(0,18){$\rho(\sigma_r)=$}
     \put(40,12){\epsfysize=25pt\epsffile{figures/12-big.eps}}
     \put(82,18){$\Longrightarrow$}
  \put(120,18){$\rho(\sigma)$}
    \put(150,18){=}
    \put(155,20){\epsffile{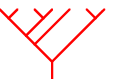}}
      \put(170,-3){\epsffile{figures/Y.eps}}
   \put(193,20){\epsffile{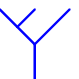}}
     \put(217,18){=}
     \put(230,7){\epsffile{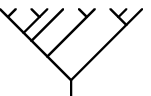}}       
\end{picture}
 \]

If $\sigma$ is 132-avoiding, then the first global descent occurs
at the maximum value of $\sigma$, and both $\sigma_l$ and $\sigma_r$ are
132-avoiding. 
Thus for 132-avoiding permutations $\sigma$, $\rho(\sigma)=\lambda(\sigma)$,
and we conclude that $\rho(\gamma(t))=t$.

\begin{lem}\label{L:rho}
 Let $n\geq 0$.
 The map $\rho\colon \frakS_n\to\calY_n$ is 
 order-preserving and for each tree $t\in\calY_n$, we have
 \begin{equation}\label{E:adjoint}
    \gamma(t)\ =\ \min\{\sigma\in\frakS_n\mid \rho(\sigma)=t\}\,.
 \end{equation}
\end{lem}
\noindent{\it Proof. }
 Suppose that $\sigma\in\frakS_n$ is not 132-avoiding.
 Then we construct a permutation $\sigma'$ with $\sigma'\lessdot\sigma$
 such that $\rho(\sigma')=\rho(\sigma)$.
 Since this process terminates when $\sigma$ is 132-avoiding, this, together
 with $\rho(\gamma(t))=t$, proves~\eqref{E:adjoint}. 

 Suppose that $\sigma$ has a 132-pattern.
 Among all 132-patterns of $\sigma$ choose a pattern 
 ($i<j<k$ with $\sigma(i)<\sigma(k)<\sigma(j)$)
 with $\sigma(k)$ maximum, and among those, with $\sigma(j)$ minimum.
 Let $m$ be the position such that $\sigma(m)=\sigma(k){+}1$.
 We must have $m\leq j$ for $m>j$ contradicts the maximality of
 $\sigma(k)$. 
 The choice of $j$ implies that either $m=j$ or else $m<i$.
 Transposing the values $\sigma(m)=\sigma(k){+}1$ and $\sigma(k)$ gives a new
 permutation  $\sigma'\lessdot\sigma$.
 We then iterate this procedure, eventually obtaining a 132-avoiding permutation.
 For example, we iterate this procedure on a permutation in $\frakS_7$:
 \[{\color{blue}
  \underline{\textcolor{black}{4}}
 \underline{\textcolor{black}{7}}
 \textcolor{black}{5} 
  \underline{\textcolor{black}{6}}
  \textcolor{black}{132}\ \to\ 
   \underline{\textcolor{black}{465}} \textcolor{black}{7132}\ \to\ 
   \textcolor{black}{4567}\underline{ \textcolor{black}{132}}\ \to\ 
    \textcolor{black}{4567213}}\,.
 \]

 We use induction on $n$ to prove that $\rho(\sigma')=\rho(\sigma)$.
 First note that $\sigma$ and $\sigma'$ have the same global descents.
 This is clear  for global descents outside of the interval
 $[m,k]$.
 By the 132-pattern at $i<j<k$, the only other possibility is if
 $\sigma$ or $\sigma'$ has a global descent at $k$, but $\sigma'$ has
 a global descent at $k$ if and only if $\sigma$ does.
 In particular, $\sigma$ and $\sigma'$ have the same first global descent.
 If this is at $k$, then $\sigma'_r=\sigma_r$ and $\sigma'_l=\sigma_l$, as only one
 of the transposed values $\sigma(k)$ and $\sigma(k)-1=\sigma(m)$ occurs to
 the left of $k$ and neither occurs to the right.
 If the first global descent is outside of the interval $[i,k]$, 
 then one of $\sigma_l$ or $\sigma_r$ contains the pattern used to construct
 $\sigma'$.  
 If $\sigma_l$ contains that pattern, then $\sigma'_r=\sigma_r$, so 
 $\rho(\sigma'_r)=\rho(\sigma_r)$, and our inductive hypothesis implies that
 $\rho(\sigma'_l)=\rho(\sigma_l)$. 
 We reach the same conclusion if $\sigma'_l=\sigma_l$, and so we
 conclude that $\rho(\sigma')=\rho(\sigma)$.

 To see that $\rho$  is order-preserving, let $\tau\leq\sigma$.
 Since $\gamma(\rho(\tau))\leq\tau$ and $\gamma(\rho(\tau))$ is 132-avoiding, we
 need only show that if $\tau\leq\sigma$ and $\tau$ is 132-avoiding, then
 $\tau\leq\gamma(\rho(\sigma))$. 
 Consider the construction of the permutation $\sigma'\lessdot\sigma$ in the
 preceding paragraph.
 Observe that $\sigma$ has exactly one more inversion, $(m,k)$, than does $\sigma'$, 
 where either $m=j$ or $m<i$. 
 Since $i<j<k$ is a 132-pattern in $\sigma$ and $\tau\leq\sigma$ is
 132-avoiding, we must have that $\tau(i)<\tau(j)<\tau(k)$.
 If $(m,k)$ were an inversion of $\tau$, then $\tau(m)>\tau(j)$,
 and  so $(m,j)$ is an inversion of $\sigma$, as $\tau\leq\sigma$ implies that
 $\Inv(\tau)\subseteq\Inv(\sigma)$.
 But this contradicts the choice of $j$, which implies that 
 $\sigma(j)\geq\sigma(k){+}1=\sigma(m)$.
 We conclude by induction that $\rho$ is order-preserving. \quad\QED\medskip

\noindent{\it Proof of Theorem~$\ref{T:galois}$. }
 The result for the horizontal maps appears in~\cite{AS02}
 (Propositions 2.11 and 2.13).%; see also~\cite{BW97}.
 
 We first treat the maps between $\frakS_n$ and $\calY_n$.
 Suppose that $\sigma\in\frakS_n$ and $t\in\calY_n$ satisfy 
 $\lambda(\sigma)\leq t$.
 Then by~\eqref{E:gamma-lambda}, we have $\sigma\leq\gamma(\lambda(w))$.
 Since $\gamma$ is order-preserving, we have $\gamma(\lambda(\sigma))\leq\gamma(t)$,
 and so $\sigma\leq\gamma(t)$.
 Conversely, suppose that $\sigma\leq\gamma(t)$.
 Since $\lambda$ is order-preserving, $\lambda(\sigma)\leq\lambda(\gamma(t))$.
 But $\lambda\circ\gamma$ is the identity, so we conclude that 
 $\lambda(\sigma)\leq t$. Thus $\lambda$ is left adjoint to $\gamma$.

 Virtually the same argument using Lemma~\ref{L:rho} shows that for 
 $t\in\calY_n$ and $\sigma\in\frakS_n$,
 \[
   \gamma(t)\leq\sigma\ \Longleftrightarrow\  t\leq\rho(\sigma)\,.
 \]
 Lastly, the remaining two equivalences, that for $t\in\calY_n$ and
 $\setS\in\calQ_n$, we have
 \begin{eqnarray*}
   L(t)\leq \setS &\Longleftrightarrow& t\leq C(\setS)\\
   C(\setS)\leq t &\Longleftrightarrow& \setS\leq R(t)
 \end{eqnarray*}
 follow from the corresponding facts for the horizontal maps via 
 $\gamma\colon\calY_n\hookrightarrow\frakS_n$. \quad\QED\medskip

We remark that the Galois connection $(\lambda,\gamma)$ can be traced back to~\cite[Sec. 9]{BW97}. Generalizations appear in~\cite{Re04a}. All the ingredients for the Galois connection $(L,C)$ also appear in~\cite[Prop. 3.5]{LR02}.

%%%%%%%%%%%%%%%%%%%%%%%%%%%%%%%%%%%%%%%%%%%%%%%%%%%%%%%%%%%%%%%%%%%%
\section{Some Hopf morphisms involving $\YSym$}\label{S:Hopfmaps}

Consider the diagram
\begin{equation}\label{E:Des-Lambda-L}
  \raisebox{-20pt}{  \begin{picture}(133,36)(-10,-1)

          \put(48,28){$\YSym$}
    \put(-10,0){$\QSym$}    \put(107,0){$\SSym$}
    \put(45,27){\vector(-2,-1){30}} \put(45,27){\vector(-2,-1){25}}
    \put(102.5, 2){\vector(-1,0){75}}\put(102.5, 2){\vector(-1,0){70}}
    \put(115,12){\vector(-2,1){30}}\put(115,12){\vector(-2,1){25}}
  
    \put(27,24){\scriptsize$\calL$}    \put(102,22){\scriptsize$\Lambda$}
          \put(61,5){\scriptsize$\calD$}
   \end{picture} }
\end{equation}
These are surjective morphisms of Hopf algebras (Proposition~\ref{P:ontoHopf}) and the diagram commutes (on the fundamental basis, this is the commutativity of the first diagram in Theorem~\ref{T:galois}). 
We use the Galois connections of Section~\ref{S:galois} to determine the
effect of these maps on the bases of monomial functions. 

Recall that  a permutation $\sigma\in\frakS_n$ has the
form $Z(\setS)$ for some $\setS\in\calQ_n$ if and only if it is $(132,213)$-avoiding (Section~\ref{S:combinatorics}). Thus~\eqref{E:MRstuff2} states that
\[
  \calD(\calM_\sigma) = \left\{\begin{array}{ccl}
             \setM_{\Des(\sigma)}&\ &\textrm{if  $\sigma$ is $(132,213)$-avoiding,} \\
                 0   &  &\textrm{otherwise.}
            \end{array}\right.
\]
We derive similar descriptions for the maps $\Lambda$ and $\calL$ in terms of
pattern avoidance. 
We say that a tree $t\in\calY_n$ is
 \epsfysize=10pt\epsffile{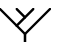}-avoiding if it 
has the form 
$t=C(\setS)$ for some $\setS\in\calQ_n$. Since $\gamma\circ C=Z$, 
the tree $t$ is  \epsfysize=10pt\epsffile{figures/213.a.eps}-avoiding if and
only if the permutation $\gamma(t)$ is $213$-avoiding. 
Geometrically, $t$ is  \epsfysize=10pt\epsffile{figures/213.a.eps}-avoiding if
every leftward pointing branch  emanates from the rightmost branch.
Equivalently, if each indecomposable component is a right comb.
For example, the two trees on the left below are
 \epsfysize=10pt\epsffile{figures/213.a.eps}-avoiding, while the tree on the
right is not. 
 \[
  \epsffile{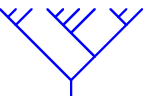}\qquad
  \epsffile{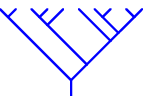}\qquad
  \epsffile{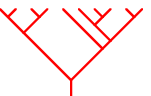}
 \]

\begin{thm}\label{T:maps}
  Let $\sigma\in\frakS_n$ and $t\in\calY_n$.
  Then
 \begin{align}
 \label{E:maps1}
  \Lambda(\calM_\sigma) &= \left\{\begin{array}{lcl}
         M_{\lambda(\sigma)}&\ &\textrm{if $\sigma$ is $132$-avoiding,}\\
            0&&\textrm{otherwise;} \end{array}\right.\\
 \label{E:maps2}
   \calL(M_t) &= \left\{\begin{array}{lcl}
         \setM_{L(t)}&\ &\textrm{if $t$ is
         \epsfysize=10pt\epsffile{figures/213.a.eps}-avoiding,}\\ 
           0&&\textrm{otherwise.}\end{array}\right.
 \end{align}
\end{thm}

\noindent{\it Proof. }
 If we have poset maps $f\colon P\to Q$ and $g\colon Q\to P$ with $f$
 left adjoint to $g$, then Rota~\cite[Theorem 1]{Rot} showed that the M\"obius 
 functions of $P$ and $Q$ are related by
 \begin{equation*}
   \forall\ x\in P \text{ and }\ w\in Q,\quad \ 
      \sumsub{ y\in P\\x\leq y,\, f(y)=w}\!\!\mu_P(x,y)\ =\ 
      \sumsub{v\in Q\\v\leq w,\, g(v)=x}\!\!\mu_Q(v,w)\, .
 \end{equation*}
 Thus, for $\sigma\in\frakS_n$ and $t\in\calY_n$, we have 
 \[
    \sumsub{\sigma\leq\tau\in\frakS_n\\\lambda(\tau)=t\rule{0pt}{10pt}}
       \mu_{\frakS_n}(\sigma,\tau)\ =\ 
    \sumsub{s \leq t\in\calY_n\\\gamma(s)=\sigma\rule{0pt}{10pt}}
               \mu_{\calY_n}(s,t)\,.
 \]
 If $\sigma$ is not 132-avoiding, then the index set on the right hand side is
 empty.  
 If $\sigma$ is 132-avoiding, then by~\eqref{E:gamma-lambda}, the index set consists
 only of the tree $s=\lambda(\sigma)$, so we have
 \begin{equation}\label{E:mobius}
   \sumsub{\sigma\leq\tau\in\frakS_n\\\lambda(\tau)=t}
          \mu_{\frakS_n}(\sigma,\tau)\ =\ 
    \begin{cases}
         \mu_{\calY_n}(\lambda(\sigma),t) & \text{if $\sigma$ is 132-avoiding,}\\
          0 & \text{if not.}
   \end{cases}
 \end{equation}
%

% We complete the proof of Theorem~\ref{T:maps}.
% Recall~\eqref{E:def-monomial} that 
%
% \[
%      \calM_\sigma\ =\ \sum_{\sigma\leq\tau} 
%         \mu_{\frakS_n}(\sigma,\tau)\cdot \calF_\tau\,,
% \]
%
% hence
%
Now, according to~\eqref{E:def-monomial} and~\eqref{E:deflambda},
 \begin{align*}
   \Lambda(\calM_\sigma)
     &\ =\ \sum_{\sigma\leq\tau}\mu_{\frakS_n}(\sigma,\tau)\cdot F_{\lambda(\tau)}
%     \\ &\ =\ 
     \quad = \quad 
     \sum_{t}\Bigl(\sumsub{\sigma\leq\tau\\\lambda(\tau)=t}
                     \mu_{\frakS_n}(\sigma,\tau) \Bigr)\cdot F_{t}\\
     &\ =\ \begin{cases}
   \sum_{t}\mu_{\calY_n}(\lambda(\sigma),t)\cdot F_t & 
          \text{if $\sigma$ is 132-avoiding}\\
          0 & \text{if not.}
  \end{cases}
\end{align*}
This proves the first part of the theorem in view of~\eqref{E:def-Ymonomial}.
%
% \[ 
%   M_{\lambda(\sigma)}\ =\ \sum_{\lambda(\sigma)\leq t} 
% \mu_{\calY_n}(\lambda(\sigma),t)\cdot F_t\,,
% \]
%
 %which proves the first part of the theorem.

 For the second part, let $t\in\calY_n$. We just showed that $\Lambda(\calM_{\gamma(t)})=M_t$. {}From the commutativity of~\eqref {E:Des-Lambda-L} we deduce
 \[
   \calL(M_t) \ =\ \calL\bigl(\Lambda(\calM_{\gamma(t)})\bigr)\ =\ 
     \calD(\calM_{\gamma(t)})\,.
 \]
 The remaining assertion follows from this and~\eqref{E:MRstuff2}. \quad\QED

%%%%%%%%%%%%%%%%%%%%%%%%%%%%%%%%%%%%%%%%%%%%%%%%%%%%%%%%%%%%%%%%%%%%%%%%%%%%
\section{Geometric interpretation of the product of $\YSym$}\label{S:product}

For $\SSym$, the  multiplicative structure constants with
respect to the basis $\{\calM_\sigma\mid\sigma\in\frakS_\infty\}$ 
were given a combinatorial
description in Theorem~4.1 of~\cite{AS02}.
In particular, they are non-negative.
An immediate consequence of this and~\eqref{E:maps1} is that the Hopf algebra
$\YSym$ has non-negative multiplicative structure constants with
respect to the basis $\{M_t\mid t\in\calY_\infty\}$.
We give a direct combinatorial interpretation of
these structure constants and complement it with a geometric
interpretation in terms of the facial structure of the associahedron.

For each permutation $\zeta\in\frakS^{(p,q)}$ consider the maps 
\[
  \varphi_\zeta\colon\frakS_p\times\frakS_q\to\frakS_{p+q}
    \qquad\text{and}\qquad
  f_\zeta\colon\calY_p\times\calY_q\to \calY_{p+q}
\]
 defined by
 \begin{equation}\label{E:fzeta}
  \varphi_\zeta(\sigma,\tau)\ :=\ (\sigma/\tau)\cdot\zeta^{-1}
   \qquad\text{and}\qquad
  f_\zeta(s,t)\ :=\ \lambda\bigl(\gamma(s)/\gamma(t)\cdot\zeta^{-1}\bigr)\,.
 \end{equation}
We supress the dependence of $\varphi_\zeta$ and $f_\zeta$ on $p$ and $q$; only
when $\zeta$ is the identity permutation does this matter.  
We view $\frakS_p\times\frakS_q$ as the Cartesian product of the posets
$\frakS_p$ and $\frakS_q$, and similarly for $\calY_p\times\calY_q$. 
The map $\varphi_\zeta$ is order-preserving~\cite[Prop. 2.7]{AS02}.
Since 
 \begin{equation}\label{E:f-varphi}
  f_\zeta(s,t)\ =\ (\lambda\circ\varphi_\zeta)(\gamma(s),\gamma(t))\,,
 \end{equation}
$f_\zeta$ is also order-preserving.

We describe the structure constants of $\SSym$ and $\YSym$ in
terms of these maps. 

\begin{prop}[Theorem~4.1 of~\cite{AS02}]\label{P:4.1}
 Let $\sigma\in\frakS_p$, $\tau\in\frakS_q$, and $\rho\in\frakS_{p+q}$.
 Then the coefficient of $\calM_\rho$ in the product 
 $\calM_\sigma\cdot\calM_\tau$ is
\[
  \alpha^\rho_{\sigma,\tau}\ :=\ 
  \#\{\zeta\in\frakS^{(p,q)}\mid (\sigma,\tau)=\max 
                \varphi_\zeta^{-1}[\id_{p+q},\rho]\}\,,
\]
where $[1_{p+q},\rho]=\{\rho'\in\calY_{p+q}\mid \rho'\leq \rho\}$, the interval in
 $\frakS_{p+q}$ below $\rho$.
\end{prop}

\begin{thm}\label{T:combprod}
 Let $s\in\calY_p$, $t\in\calY_q$, and $r\in\calY_{p+q}$.
 The coefficient of $M_r$ in $M_s\cdot M_t$ is
 \begin{equation}\label{E:combprod}
   \#\{\zeta\in\frakS^{(p,q)}\mid (s,t)=\max f_\zeta^{-1}[1_{p+q},r]\}\,.
 \end{equation}
 \end{thm}

\noindent{\it Proof. }
 By~\eqref{E:maps1},
 $M_s\cdot M_t=\Lambda\bigl(\calM_{\gamma(s)}\cdot\calM_{\gamma(t)}\bigr)$.
 We evaluate this using Proposition~\ref{P:4.1} and~\eqref{E:maps1}
 to obtain
\[
  M_s\cdot M_t\ =\ \Lambda\Bigl(\sum_{\rho\in\frakS_{p+q}} 
        \alpha^\rho_{\gamma(s),\gamma(t)}\,\calM_\rho \Bigr) \ =\ 
  \sum_{r\in\calY_{p+q}} \alpha^{\gamma(r)}_{\gamma(s),\gamma(t)} M_r\,.
\]
 The constant $\alpha^{\gamma(r)}_{\gamma(s),\gamma(t)}$ is equal to
\[
  \#\{\zeta\in\frakS^{(p,q)}\mid (\gamma(s),\gamma(t))=\max 
                \varphi_\zeta^{-1}[\id_{p+q},\gamma(r)]\}\,.
\]
 The Galois connections between $\frakS_{p+q}$ and $\calY_{p+q}$ of
 Theorem~\ref{T:galois} imply that 
\[
  \lambda\bigl(\gamma(s)/\gamma(t)\cdot\zeta^{-1}\bigr)\ \leq\ r
   \ \Longleftrightarrow\ 
    \gamma(s)/\gamma(t)\cdot\zeta^{-1}\ \leq\ \gamma(r)\,.
\]
By the definitions of $f_\zeta$ and $\varphi_\zeta$, it follows that
\[
  (s,t)\ =\ \max f_\zeta^{-1}[1_{p+q},r]\ \iff\ 
  (\gamma(s),\gamma(t))\ =\ \max \varphi_\zeta^{-1}[\id_{p+q},\gamma(r)]\,,
\]
 and hence $\alpha^{\gamma(r)}_{\gamma(s),\gamma(t)}$ equals~\eqref{E:combprod},
 which completes the proof.
\quad\QED\medskip

\begin{figure}[htb]
 \[
   \epsffile{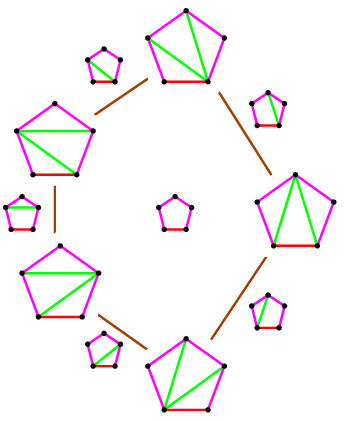}\qquad\quad
   \epsffile{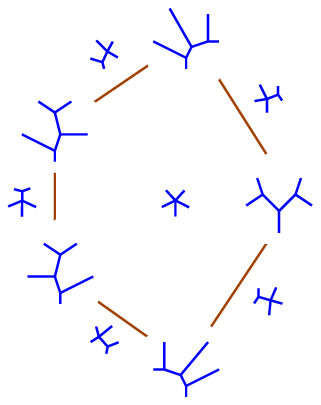}
%sottile:figures: fig2dev -Leps -m0.35 A3.tree.fig  A3.tree.eps
%sottile:figures: fig2dev -Leps -m0.35 A3.tri.fig  A3.tri.eps
 \]
\caption{Two views of the associahedron $\calA_3$}
\label{F:assoc}
\end{figure}

The Hasse diagram of $\calY_n$ is isomorphic to the 1-skeleton of 
the {\it associahedron} $\calA_n$, an $(n{-}1)$-dimensional polytope.
(See~\cite[pp.~304--310]{Zi95} and~\cite[p.~271]{St99}.)
The faces of $\calA_n$ are in one-to one correspondence with
collections of non-intersecting diagonals of a polygon with $n{+}2$ sides (an
($n{+}2$)-gon). 
Equivalently, the faces of $\calA_n$ correspond to  polygonal subdivisions of an
$n{+}2$-gon with  
facets corresponding to diagonals and vertices to triangulations.
The dual graph of a polygonal subdivision is a planar tree and the dual graph
of a triangulation is a planar binary tree.
If we distinguish one edge to be the root edge, the trees are rooted, and this
furnishes a bijection between the vertices of $\calA_n$ and $\calY_n$.
Figure~\ref{F:assoc} shows two views of the associahedron $\calA_3$,
the first as polygonal subdivisions of the pentagon, and the second as 
the corresponding dual graphs (planar trees).
The root is at the bottom.
%%%%%%%%%%%%%%%%%%%%%%%%%%%%%%%%%%%%%%%%%%%%%%%%%%%%%%%%%%%%%%%%%%%%%%%%
%%%%%%%%%%%%%%%%%%%%%%%%%%%%%%%%%%%%%%%%%%%%%%%%%%%%%%%%%%%%%%%%%%%%%%%%

We describe the map $\lambda\colon\frakS_n\to\calY_n$ in terms of triangulations
of the ($n{+}2$)-gon where we label the vertices with $0,1,\dots,n,n{+}1$  
beginning with the left vertex of the root edge and proceeding clockwise.
Let $\sigma\in\frakS_n$ and set 
$w_i:=\sigma^{-1}(n+1-i)$, for $i=1,\ldots,n$.
%$w:=(w_1,w_2,\dotsc,w_n)=\sigma^{-1}(n,n{-}1,\dotsc,2,1)$.
This records the positions of the values of $\sigma$ taken in decreasing
order.
We inductively construct the triangulation, beginning with the empty
triangulation consisting of the root edge, and after $i$ steps we have a
triangulation $T_i$ of the polygon
\[
   P_i\ :=\ \textit{Conv}\{0,n{+}1,w_1,\dotsc,w_i\}\,.
\]
Some edges of $P_i$ will be edges of the original ($n{+}2$)-gon and others
will be diagonals.
Each diagonal cuts the ($n{+}2$)-gon into two pieces, one containing
$P_i$ and the other a polygon which is not yet triangulated and
whose root edge we take to be that diagonal. 
Subsequent steps add to the triangulation $T_i$ and its support $P_i$.

First set $T_1:=\textit{Conv}\{0,n{+}1,w_1\}$, the triangle with base the root
edge and apex the vertex $w_1=\sigma^{-1}(n)$.
Set $P_1:=T_1$ and continue.
After $i$ steps we have constructed $T_i$ and $P_i$ in such a way that the vertex $w_{i+1}$  is
not in $P_i$. Hence it must lie in some untriangulated polygon consisting of some consecutive edges of the ($n{+}2$)-gon and a diagonal that is an edge of
$P_i$. 
Add the join of the vertex $w_{i+1}$ and the diagonal to the triangulation to
obtain a triangulation $T_{i+1}$ of the polygon $P_{i+1}$.
The process terminates when $i=n$.

For example, we display this process for the permutation $\sigma=316524$,
where we label the vertices of the first octagon:
\[
  \begin{picture}(400,60)
   \put(0,0){
    \begin{picture}(70,60)
    \put(8,0){\epsffile{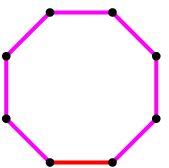}}
%fig2dev -Lps -m0.4 1256.fig 1256.eps 
    \put(18,50){3}   \put(38,50){4}
    \put( 0,29){2}   \put(58,29){5}
    \put( 1,10){1}   \put(58,10){6}
   \end{picture}
}
  \put(75,21){\large$\longmapsto$}
  \put(90,0){
   \begin{picture}(70,60)
    \put(8,0){\epsffile{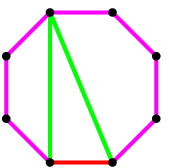}}
   \end{picture}
}
  \put(155,21){\large$\longmapsto$}
  \put(170,0){
   \begin{picture}(70,60)
    \put(8,0){\epsffile{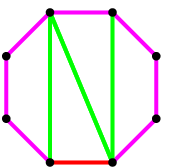}}
   \end{picture}
}
  \put(235,21){\large$\longmapsto$}
  \put(250,0){
   \begin{picture}(70,60)
    \put(8,0){\epsffile{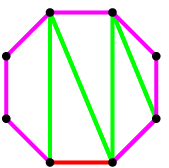}}
   \end{picture}
}
  \put(315,21){\large$\longmapsto$}
  \put(330,0){
   \begin{picture}(70,60)
    \put(8,0){\epsffile{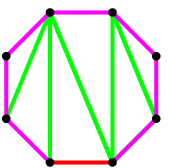}}
   \end{picture}
}
 \end{picture}
\]
The last two steps are supressed as they add no new diagonals.
The dual graph to the triangulation $T_n$ is the planar binary tree
$\lambda(\sigma)$. 
Here is the triangulation, its dual graph, and a `straightened' version, which
we recognize as the tree $\lambda(316524)$.
\[
   \epsffile{figures/T4.eps}\qquad
    \epsffile{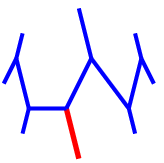}\qquad
    \epsffile{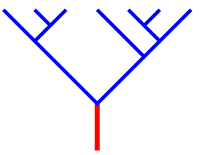}
    %fig2dev -Leps -m0.5 T4-straight.fig T4-straight.eps
\]

A subset $\setS$ of $[n]$ determines a face  $\Phi_\setS$ of the associahedron
$\calA_n$ as follows.
Suppose that we label the vertices of the ($n{+}2$)-gon as above.
Then the vertices labeled $0,n{+}1$ and those labeled by $\setS$ form a 
($\#\setS+2$)-gon whose edges include a set $E$ of non-crossing diagonals of
the original ($n{+}2$)-gon.
These diagonals determine the face $\Phi_\setS$ of $\calA_n$ corresponding to
$\setS$. 
We give two examples of this association when $n=6$ below.
\[
  \begin{picture}(350,60)
   \put(0,20){$\{1,2,5,6\}\ \longleftrightarrow$ }
   \put(84,0){
    \begin{picture}(70,60)
    \put(8,0){\epsffile{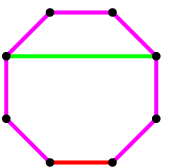}}
%fig2dev -Leps -m0.4 1256.fig 1256.eps 
    \put(18,50){3}   \put(38,50){4}
    \put( 0,29){2}   \put(58,29){5}
    \put( 1,10){1}   \put(58,10){6}
   \end{picture}
}
  \put(200,20){$\{2,4,5\}\ \longleftrightarrow$ }
  \put(274,0){
   \begin{picture}(75,72)
    \put(9,0){\epsffile{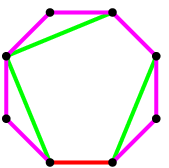}}
%fig2dev -Leps -m0.4 245.fig  245.eps  
    \put(18,50){3}   \put(38,50){4}
    \put( 0,29){2}   \put(58,29){5}
    \put( 1,10){1}   \put(58,10){6}
   \end{picture}
}
%sottile:figures: fig2dev -Leps -m0.5 36.fig  36.fig   
   \end{picture}
\]
We determine the image of $f_\zeta$ using the 
above description of the map $\lambda\colon\frakS_n\to\calY_n$.

\begin{prop}\label{P:lambda_sigma}
 Given $\zeta\in\frakS^{(p,q)}$, the image of $\lambda\circ \varphi_\zeta$ coincides with the
 image of $f_\zeta$ and equals the face 
 $\Phi_{\zeta\{p{+}1,\dotsc,p{+}q\}}$ of $\calA_{p+q}$.
 This is a facet if and only if $\zeta$ is $132$-avoiding.
\end{prop}

\noindent{\it Proof. }
 Let $\sigma\in\frakS_p$ and $\tau\in\frakS_q$, and set
 $\rho:=\varphi_\zeta(\sigma,\tau)=(\sigma/\tau)\cdot\zeta^{-1}$. 
 Since the $q$ largest values of $\rho$ lie in the positions
 $\setS:=\zeta\{p{+}1,\dotsc,p{+}q\}$, the triangulation
 $\lambda(\rho)$ is obtained by
 triangulating
 $Q:=\textit{Conv}\{0,p{+}q{+}1,\setS\}$ 
 with $\lambda(\tau)$, and then placing triangulations given by parts of
 $\sigma$ in the polygons that lie outside $Q$.
 More precisely, suppose that 
 $\{1,\dotsc,p{+}q\}-\setS$
 consists of strings of $a_1, a_2, \dotsc, a_r$ consecutive numbers.
 Then the $i$th polygon $P_i$ outside of $Q$ is triangulated according to 
 $\lambda(\st(\sigma_{A_i+1},\dotsc,\sigma_{A_i+a_i}))$, 
 where $A_i=a_1+\dotsb+a_{i-1}$.
 Thus all triangulations of the $(n{+}2)$-gon that include the edges of $Q$ are
 obtained from permutations of the form 
 $(\gamma(s_1)/\gamma(s_2)/\dotsb/\gamma(s_r), \ \gamma(t))$, where 
 $s_i\in\calY_{a_i}$ and $t\in\calY_q$.
 But this describes the face 
 $\Phi_{\setS}$ and shows that 
 $f_\zeta(\calY_p\times\calY_q)
  =(\lambda\circ\varphi_\zeta)(\frakS_p\times\frakS_q)=\Phi_{\setS}$.

 This face will be a facet only if $r=1$, that is if
 $\{1,\dotsc,p{+}q\}-\setS$
 consists of consecutive numbers, which is equivalent to $\zeta$ avoiding the
 pattern 132.
\quad\QED\medskip

We say that a face of $\calA_{p+q}$ of the form $\Phi_\setS$ with $\#\setS=q$ has  {\it type} $(p,q)$. If a face has a type, this type is unique.
A permutation $\zeta\in\frakS^{(p,q)}$ is uniquely determined by the set $\zeta\{p{+}1,\dotsc,p{+}q\}$. 
Therefore, a face of type $(p,q)$ is the image of $f_\zeta$ for a unique
permutation $\zeta\in\frakS^{(p,q)}$. This allows us to speak of the vertex of the face corresponding to a pair $(s,t)\in\calY_p\times\calY_q$ (under $f_\zeta$).

We deduce a geometric interpretation of the multiplicative structure
constants from Proposition~\ref{P:lambda_sigma} and Theorem~\ref{T:combprod}.

\begin{cor}\label{C:combprodY}
 Let $s\in\calY_p$ and $t\in\calY_q$.
 The coefficient of $M_r$ in the product $M_s\cdot M_t$ equals the number of
 faces of the associahedron $\calA_{p+q}$ of type $(p,q)$ such 
 %with the property
 that the vertex corresponding to $(s,t)$ is below $r$, and it is the
 maximum vertex on its face below $r$.
\end{cor}

Proposition~\ref{P:lambda_sigma} describes the image of a
facet $\varphi_\zeta(\frakS_p\times\frakS_q)$ of the permutahedron for
$\frakS_{p+q}$ under the map $\lambda$. 
More generally, it is known  that the image of any face of the
permutahedron is a face of the associahedron~Tonks~\cite{To97}.
The map from the permutahedron to the associahedron can also be
understood by means of the theory of {\em fiber polytopes}~\cite[Sec. 5]{BS94},
~\cite[Sec. 4.3]{Rei02}. 
For more on the permutahedron and associahedron, see~\cite{Lod04}.

%%%%%%%%%%%%%%%%%%%%%%%%%%%%%%%%%%%%%%%%%%%%%%%%%%%%%%%%%%%%%%%%%%%%%
\section{Cofreeness and the coalgebra structure of $\YSym$}

We compute the coproduct on the basis $\{M_t\mid y\in\calY_\infty\}$ and deduce
the existence of a new grading  for which $\YSym$ is cofree. 
We show that $\{M_{t\vee |}\mid t\in\calY_\infty\}$ is a basis for the space of
 primitive elements and describe the coradical filtration of $\YSym$.
Since $\YSym$ is the graded dual of the Loday-Ronco Hopf algebra $\textit{LR}$,
 this work strengthens the result~\cite[Theorem 3.8]{LR98} of Loday and Ronco
 that $\textit{LR}$ is a free associative algebra. 

\begin{thm}\label{T:coproduct}
 Let $r\in\calY_n$.
 Then
\begin{equation}\label{E:coproduct}
   \Delta(M_r)\ =\ \sum_{r= s\backslash t} M_s\ten M_t\,.
\end{equation}
\end{thm}

\noindent{\it Proof. }
 Suppose that $\rho\in\frakS_n$ is 132-avoiding and we decompose $\rho$ as 
 $\rho=\sigma\backslash\tau$.
 Then both $\sigma$ and $\tau$ are 132-avoiding: a 132 pattern in either
 would give a 132 pattern in $\rho$.
 We use that $\Lambda$ is a morphism of coalgebras and~\eqref{E:MRstuff1} to obtain
 \begin{eqnarray*}
  \Delta(M_r)&=& \Delta(\Lambda(\calM_{\gamma(r)}))\ =\ 
                  \Lambda(\Delta(\calM_{\gamma(r)}))\\
        &=& \Lambda \Bigl(\sum_{\gamma(r)=\sigma\backslash\tau}
                  \calM_\sigma\ten\calM_\tau\Bigr)\\
        &=& \sum_{\gamma(r)=\sigma\backslash\tau} 
                  M_{\lambda(\sigma)}\ten M_{\lambda(\tau)}
        \quad =\quad  \sum_{r=s\backslash t} M_s\ten M_t\,,
 \end{eqnarray*}
 the last equality by Proposition~\ref{P:starslash}. \quad\QED\medskip 

We recall the notion of cofree graded coalgebras. 
Let $V$ be a vector space and set
 \[ 
   Q(V)\ :=\ \bigoplus_{k\geq 0} V^{\ten k}\,,
 \]
which is naturally graded by $k$.
Given $v_1,\dotsc,v_k\in V$, let $v_1\iten v_2\iten\dotsb\iten v_k$ denote the corresponding tensor in $V^{\ten k}$. 
 Under the \emph{deconcatenation} coproduct
 \[
   \Delta(v_1\iten\dotsb\iten v_k)\ =\ 
   \sum_{i=0}^k (v_1\iten\dotsb\iten v_i)\otimes(v_{i+1}\iten\dotsb\iten v_k)\,,
 \]
and counit $\epsilon(v_1\iten\dotsb\iten v_k)=0$ for $k\geq 1$, 
$Q(V)$ is a graded connected coalgebra, the {\em cofree graded coalgebra on $V$}.

We show that $\YSym$ is cofree by first defining
a second coalgebra grading, where the degree of $M_t$ 
is the number of branches of $t$ emanating from the rightmost branch
(including the leftmost branch), that is, $1+\#R(t)$.
This is also the number of components in the (right) decomposition of $t$ into
progressive trees $t=t_1\backslash t_2 \backslash \cdots \backslash t_k$, as in
Section~\ref{S:basic}.

First, set $\calY^{0}:=\calY_0$, and for
$k\geq 1$, let
 \begin{align*}
   \calY_n^{k}&\ :=\ \{t\in\calY_n\mid   t
                \text{ has exactly $k$ progressive components}\}, \
       \text{ and }\\
    \calY^{k}&\ :=\ \coprod_{n\geq k}\calY_n^{k}\,.
 \end{align*}
In particular $\calY^1$ consists of the progressive trees, those of the form
$t\vee |$.
For instance,
\begin{eqnarray*}
  \calY^1&=& \bigl\{\epsffile{figures/1.eps}\bigr\}\: \cup\: 
                \bigl\{\epsffile{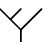}\bigr\}\: \cup\: 
                \bigl\{\epsffile{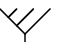},
                      \epsffile{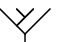}\bigr\}\: \cup\: 
               \bigl\{\epsffile{figures/1234.eps},
                      \epsffile{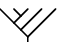},
                      \epsffile{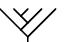},
                      \epsffile{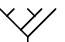},
                      \epsffile{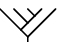}\bigr\}\: \cup\: \dotsb\\
   \calY^2&=&\bigl\{\epsffile{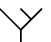}\bigr\}\: \cup\: 
                \bigl\{\epsffile{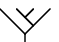},
                      \epsffile{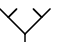}\bigr\}\: \cup\: 
               \bigl\{\epsffile{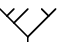},
                      \epsffile{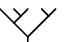},
                      \epsffile{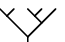},
                      \epsffile{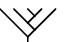},
                      \epsffile{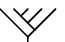}\bigr\}\: \cup\: \dotsb\\
   \calY^3&=&\bigl\{\epsffile{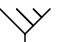}\bigr\}\: \cup\: 
                \bigl\{\epsffile{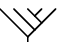},\!
                      \epsffile{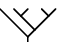},\!
                      \epsffile{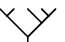}\bigr\}\: \cup\:  
               \bigl\{\epsffile{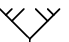},\!
                      \epsffile{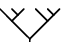},\!
                      \epsffile{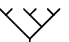},\!
                      \epsffile{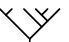},\!
                      \epsffile{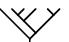},\!
                      \epsffile{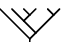},\!
                      \epsffile{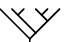},\!
                      \epsffile{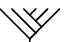},\!
                      \epsffile{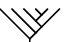}\bigr\}\: \cup\:
                \dotsb\smallskip 
\end{eqnarray*}

Let $(\YSym)^k$ be the vector subspace of $\YSym$ spanned by
$\{M_t\mid t\in\calY^k\}$.

\begin{thm}\label{T:cofree} 
 The decomposition $\YSym=\oplus_{k\geq 0}(\YSym)^k$ is a coalgebra grading. 
 With this grading $\YSym$ is a cofree graded coalgebra.
 \end{thm}

\noindent{\it Proof. }
 Let $V:=(\YSym)^1$, the span of $\{M_t\mid t\ \textrm{is progressive}\}$.
 Then the map
\[
    M_{t_1}\ten M_{t_2}\ten \dotsb \ten M_{t_k}\ \longmapsto\ 
       M_{t_1\backslash t_2\backslash\dotsb\backslash t_k}\,,
\]
 identifies $V^{\otimes k}$ with $(\YSym)^k$.
 Together with the coproduct formula~\eqref{E:coproduct}, 
 this identifies $\YSym$ with the deconcatenation coalgebra $Q(V)$.
\quad\QED\medskip

The coradical $C^0$ of a graded connected coalgebra $C$ is the
1-dimensional component in degree 0. 
The primitive elements of $C$ are
 \[
   \text{P}(C)\ :=\ \{x\in C\mid \Delta(x)=x\ten 1+1\ten x\}\,.
 \]
Set $C^1:=C^0\oplus \text{P}(C)$, the first level of the
coradical filtration.
More generally, the $k$-th level of the coradical filtration is
 \[
   C^k\ :=\ \bigl(\Delta^k\bigr)^{-1}
        \Bigl(\sum_{i+j=k}C^{\ten i}\ten C^0\ten C^{\ten j}\Bigl)\,.
 \]
We have $C^0\subseteq C^1\subseteq C^2
          \subseteq\dotsb\subseteq C=\bigcup_{k\geq 0}C^k$,
and
 \[
   \Delta(C^k)\ \subseteq\ \sum_{i+j=k}C^i\ten C^j\,.
 \]
Thus, the coradical filtration measures the complexity of
iterated coproducts.

When $C$ is a cofree graded coalgebra $Q(V)$, its
space of primitive elements is just $V$, and the $k$-th level
of its coradical filtration is $\oplus_{i=0}^k V^{\ten i}$. 
We record these facts for $\YSym$.

\begin{cor}\label{C:coradical}
 A linear basis for the $k$-th level of the coradical filtration of $\YSym$ is
 \[
   \{M_t\mid t\in \calY^k\}\,.
 \]
 In particular, a linear basis for the space of primitive elements is
 \[
    \{ M_t\mid t\ \textrm{is progressive}\}.
 \]
\end{cor}

\begin{rem}\label{R:primitive} Recall that a tree $t=t_l\vee t_r$  is progressive if and only if $t_r=1_0=|$ (Section~\ref{S:com-per-tre}). It follows that the number of progressive trees in $\calY_n$ is $\dim(\YSym)_n^1=c_{n-1}$. Theorem~\ref{T:cofree}
implies that the Hilbert series of $\YSym$ and $(\YSym)^1$ are related by
\begin{equation}\label{E:Catalan}
\sum_{n\geq 0} c_n t^n= \frac{1}{1-\sum_{n\geq 1} c_{n-1} t^n}\,.
\end{equation}
This is equivalent to the usual recursion for the Catalan numbers
\[c_n=\sum_{k=0}^{n-1} c_kc_{n-1-k}\quad \forall\, n\geq 1\,, \quad c_0=1\,.\]
\end{rem}

\begin{rem}\label{R:GL}
Let $\gr\YSym$ be the associated graded Hopf algebra to $\YSym$ under its
coradical filtration.
This is bigraded, as it also retains the original grading of $\YSym$.
Greg Warrington showed that this is commutative, and it is in fact 
the shuffle Hopf algebra generated by the $M_t$ for $t$ a progressive tree.
Grossman and Larson~\cite{GL89} defined a graded cocommutative Hopf algebra of
planar trees, whose dual is isomorphic to $\gr\YSym$~\cite{AS04}.
\end{rem}
%%%%%%%%%%%%%%%%%%%%%%%%%%%%%%%%%%%%%%%%%%%%%%%%%%%%%%%%%%%%%%%%%%%%%
\section{Antipode of $\YSym$}\label{S:antipode}

We give an explicit formula for the antipode of $\YSym$.
This is a simple consequence of the formula in~\cite[Thm. 5.5]{AS02} for the
antipode of $\SSym$.
(See Remark~9.5 in~\cite{Re04b}.)
Let $\tau\in\frakS_n$. Subsets $\setR\subseteq\GDes(\tau)$ correspond to decompositions of $\tau$ 
\[
   \setR\ \leftrightarrow\ \tau\ =\ 
    \tau_1\backslash \tau_2\backslash \dotsb \backslash \tau_r\,.
\]
For such a partial decomposition $\setR$, set
 \[
   \tau_{\setR}\ :=\ \tau_1/\tau_2/\dotsb/\tau_r\,.
 \]
For example, for the decomposition $\setR$ of the permutation
 \[
 \tau\ =\ 798563421\ =\ 132\backslash 3412\backslash 21
 \]
we have
 \[
\tau_\setR\ =\  132/ 3412/ 21\ =\ 132674598\,.
 \]
Lastly, given $n$ and a subset $\setS$ of $[n-1]$, $\frakS^\setS$ denotes the
set of permutations $\sigma\in \frakS_n$ such that $\Des(\sigma)\subseteq\setS$;
equivalently, $\frakS^\setS=[\id_n,Z(\setS)]$. 

\begin{thm}\label{T:antipodeY}
 For $t\in\calY_n$, 
\[
  S(M_t)\ =\ -(-1)^{\#R(t)}\sum_{s\in\calY_n} \kappa(t,s)M_s\ ,
\]
 where $\kappa(t,s)$ records the number of permutations $\zeta\in\frakS^{R(t)}$
 that satisfy
\begin{enumerate}
 \item[(i)]
           $\lambda(\gamma(t)_{R(t)}\cdot\zeta^{-1})\leq s$,
 \item[(ii)]
           $t\leq t'$ and $\lambda(\gamma(t')_{R(t)}\cdot\zeta^{-1})\leq s$
           implies that $t=t'$, and  
 \item[(iii)]
           If\/ $\Des(\zeta)\subseteq\setR\subseteq R(t)$ and 
           $\lambda(\gamma(t)_\setR\cdot\zeta^{-1})\leq s$, then $\setR=R(t)$.
\end{enumerate}
\end{thm}

\noindent{\it Proof. }
  Theorem~5.5 of~\cite{AS02} gives the following formula for the antipode
 on $\SSym$:
\[
  S(\calM_\tau)\ =\ -(-1)^{\#\GDes(\tau)}
   \sum_{\sigma\in\frakS_n} k(\tau,\sigma)\calM_\sigma\,,
\]
 where $k(\tau,\sigma)$ records the number of permutations 
 $\zeta\in\frakS^{\GDes(\tau)}$ that satisfy
\begin{enumerate}
 \item[(a)]
           $\tau_{\GDes(\tau)}\cdot\zeta^{-1}\leq\sigma$,
 \item[(b)]
           $\tau\leq\tau'$ and $\tau'_{\GDes(\tau)}\cdot\zeta^{-1}\leq\sigma$ 
           implies that $\tau=\tau'$, and 
 \item[(c)]
           If $\Des(\zeta)\subseteq\setR\subseteq\GDes(\tau)$ and 
           $\tau_\setR\cdot\zeta^{-1}\leq\sigma$, then $\setR=\GDes(\tau)$.
\end{enumerate}
 
 Since $\Lambda\colon\SSym\to\YSym$ is a morphism of Hopf algebras, 
 $\Lambda(S(\calM_\tau))=S(\Lambda(\calM_\tau))$.
 Also,~\eqref{E:maps1} says that $\Lambda(\calM_\sigma)=0$ unless
 $\sigma=\gamma(s)$ for some $s\in\calY_n$.
 Since we also have $R(s)=\GDes(\gamma(s))$ for $s\in\calY_n$, the theorem will
 follow from this result for $\SSym$ if the set of permutations
 $\zeta\in\frakS^{R(t)}$ satisfying conditions (i), (ii), and (iii) for
 trees $s,t\in\calY_n$ in the
 statement of the theorem equals the set satisfying (a), (b), and (c) for
 $\tau=\gamma(t)$ and $\sigma=\gamma(s)$.

 Suppose that  $\Des(\zeta)\subseteq\setR\subseteq R(t)\ (=\GDes(\tau))$.
 Then 
 $\gamma(t)_{R(t)}\cdot\zeta^{-1} = \tau_{\GDes(\tau)}\cdot\zeta^{-1}$, and so
\[
  \lambda(\gamma(t)_{R(t)}\cdot\zeta^{-1})\leq s\ 
  \Longleftrightarrow\ 
  \tau_{\GDes(\tau)}\cdot\zeta^{-1}\leq\sigma\,,
\]
 by the Galois connection (Theorem~\ref{T:galois}).
 This shows that (i) and (a) are equivalent, as well as (iii) and (c).

We show that (b) implies (ii). Suppose that $t'\in\calY_n$ satisfies
\[
   t\leq t' 
     \qquad \text{and}\qquad 
   \lambda(\gamma(t')_{R(t)}\cdot\zeta^{-1})\leq s\,.
\]
Let $\tau':=\gamma(t')$. Applying $\gamma$ to the first inequality and treating
the second as in the  preceding paragraph (replacing $t'$ for $t$ ) we obtain 
\[
   \tau\leq\tau' \qquad \text{and}\qquad \tau'_{\GDes(\tau)}\cdot\zeta^{-1}\leq
   \sigma\,.
\]
Hypothesis (b) implies $\tau=\tau'$ and applying $\gamma$ we conclude $t=t'$, so
(ii) holds. 

 To see that (ii) implies (b), suppose that $\tau'\in\frakS_n$
 satisfies
 \begin{equation}\label{Eq:necessary}
  \gamma(t)\ \leq\ \tau'\qquad\textrm{and}\qquad
  \tau'_{R(t)}\cdot\zeta^{-1}\ \leq\ \gamma(s)\,.
 \end{equation}
 By Remark~\ref{rem:fix} below,
 $\lambda(\tau'_{R(t)}\cdot\zeta^{-1})=
 \lambda\Bigl(\gamma\bigl(\lambda(\tau')\bigr)_{R(t)}\cdot\zeta^{-1}\Bigr)$.
 If we apply $\lambda$ to~\eqref{Eq:necessary}, we obtain
\[
  t\ \leq\ \lambda(\tau')\qquad\textrm{and}\qquad
  \lambda\Bigl(\gamma\bigl(\lambda(\tau')\bigr)_{R(t)}\cdot\zeta^{-1}\Bigr)\ \leq\ s\,.
\]
Assuming (ii), we conclude that $t=\lambda(\tau')$. 
This implies $\gamma(t)\geq\tau'$ by~\eqref{E:gamma-lambda}.
Hence $\tau=\gamma(t)=\tau'$, and (b) holds. 
\quad\QED\medskip 

We remark that by similar techniques one may derive an explicit formula for
the antipode of $\YSym$ on the fundamental basis $F_t$, working from
the corresponding result  for $\SSym$~\cite[Thm.~5.4]{AS02}.

For $\tau\in\frakS_n$, let $\overline{\tau}:=\gamma(\lambda(\tau))$, the 
unique 132-avoiding permutation such that
$\lambda(\tau)=\lambda(\overline{\tau})$. 
Suppose that $1\leq a<b\leq n$.
We define a premutation $\overline{\tau}^{[a,b]}$ which has no 132-patterns
having values in the interval $[a,b]$
(no occurrences of  $i<j<k$ with
 $a\leq \tau(i)<\tau(k)<\tau(j)\leq b$), and which satisfies 
$\lambda(\tau)=\lambda(\overline{\tau}^{[a,b]})$. 
Set $\setS=\{s_1<\dotsb<s_m\}:=\tau^{-1}([a,b])$ and let $\sigma$ be the
permutation $\st(\tau(s_1),\dotsc,\tau(s_m))$, the standard permutation formed
by the values of $\tau$ in the interval $[a,b]$.
Define $\overline{\tau}^{[a,b]}\in\frakS_n$ to be the permutation 
\[
   \overline{\tau}^{[a,b]}(i)\ =\ \left\{
    \begin{array}{rcl}  \tau(i)&\quad&\mbox{\rm if } i\not\in\setS\\
      a-1+\overline{\sigma}(j)&& \mbox{\rm if }i=s_j\in\setS
    \end{array}\right.
\]

\begin{lem}\label{lem:fix}
 With the above definitions, $\lambda(\tau)=\lambda(\overline{\tau}^{[a,b]})$.
\end{lem}

\begin{rem}\label{rem:fix}
 Suppose that $\tau$ is a permutation and $\setR=\{r_1<\dotsb<r_{m-1}\}$ is a subset of $\GDes(\tau)$.
 Thus $\tau=\tau_1\backslash\tau_2\backslash\dotsb\backslash\tau_m$ with 
 $\tau_i\in\frakS_{r_i-r_{i-1}}$, where $0=r_0$ and $r_m=n$.
 Then, by Proposition~\ref{P:starslash}, 
 $\overline{\tau}=\overline{\tau_1}\backslash\dotsb\backslash\overline{\tau_m}$.
 
 Let $\zeta\in\frakS^\setR$ and consider the permutations
 \begin{equation*}
  \tau_\setR\cdot\zeta^{-1}\ =\ (\tau_1/\tau_2/\dotsb/\tau_m)\cdot\zeta^{-1}
  \quad\text{and}\quad
  \overline{\tau}_\setR\cdot\zeta^{-1}\ =\ 
      (\overline{\tau_1}/\overline{\tau_2}/\dotsb/\overline{\tau_m})\cdot\zeta^{-1}\,.
 \end{equation*}
 Observe that by the definitions preceeding the statement of the lemma,
\[
  (\tau_1/\dotsb/\overline{\tau_i}/\dotsb/\tau_m)\cdot\zeta^{-1}
  \ =\ \overline{\tau_\setR\cdot\zeta^{-1}}^{[1+r_{i-1},r_i]}\,.
\]
 Thus 
\[
   \overline{\tau}_\setR\cdot\zeta^{-1}\ =\ 
  \overline{(\dotsb(\overline{\tau_\setR\cdot\zeta^{-1}}^{[1,r_1]})\dotsb)}^{[1+r_{m-1},n]}\
  . 
\]
 By Lemma~\ref{lem:fix} we conclude that 
\[
   \lambda(\tau_\setR\cdot\zeta^{-1})\ =\ 
   \lambda(\overline{\tau}_\setR\cdot\zeta^{-1})\ =\ 
   \lambda\Bigl(\gamma\bigl(\lambda(\tau)\bigr)\Bigr)_\setR\cdot\zeta^{-1})\,,
\]
 which was needed in the proof of 
Theorem~\ref{T:antipodeY}.
\end{rem}

\noindent{\it Proof of Lemma~$\ref{lem:fix}$. }
 We prove this by  increasing induction on $n$ and decreasing induction on the
 length of the permutation $\tau$. 
 The initial cases are trivial and immediate.
 Consider $132$-patterns in $\tau$ with values in $[a,b]$.
 If $\tau$ has no $132$-pattern with values in $[a,b]$, then
 $\overline{\tau}^{[a,b]}=\tau$, and there is nothing to show.

 Otherwise, consider the $132$-patterns in $\tau$ with values in $[a,b]$ where
 $\tau(j)$ is minimal, and among those, consider patterns where $\tau(k)$ is also
 minimal.
 Finally, among those, consider the one with $\tau(i)$ maximal.
 That is $(\tau(j),\tau(k),-\tau(i))$ is minimal in the lexicographic order.
 We claim that $\tau(i)=\tau(k)-1$.
 Indeed, define $m$ by $\tau(m)=\tau(k)-1$.
 We cannot have $j<m$, for then $i<j<m$ would give a $132$-pattern with
 values in $[a,b]$ where $(\tau(j),\tau(m))$ preceeds $(\tau(j),\tau(k))$ in
 the lexicographic order.
 Since $m<j$, the choice of $i$ forces $m=i$.

 Transposing the values $\tau(i)$ and $\tau(k)$ gives a permutation $\tau'$ with $\tau\lessdot\tau'$.
 We claim that $\lambda(\tau')=\lambda(\tau)$.
 This will complete the proof, as we are proceeding by downwward induction on
 the length of $\tau$.

 We prove this by induction on $n$.
 Consider forming the trees $\lambda(\tau)$  and
 $\lambda(\tau')$. 
 Let $m:=\tau^{-1}(n)$, then $\tau'(m)=n$, also.
 As in the definition of $\lambda$~\eqref{E:def-lambda}, form $\tau_l$ and
 $\tau_r$, and the same for $\tau'$. 
 If $i<m<k$, then $\tau'_l=\tau_l$ and $\tau'_r=\tau_r$, and so 
 $\lambda(\tau)=\lambda(\tau')$.
 If $k<m$, then $\tau'_r=\tau_r$, but $\tau'_l\neq\tau_l$.
  However, $\tau'_l$ is obtained from $\tau_l$ by interchanging the values $\tau_l(i)$ and $\tau_l(k)$, and $i<j<k$ is a 132-pattern in 
 $\tau_l$ with values in an interval where 
 $(\tau_l(j),\tau_l(k),-\tau_l(i))$ is minimal in the lexicographic order. 
 By induction on $n$, $\lambda(\tau_l)=\lambda(\tau'_l)$, and so 
 $\lambda(\tau)=\lambda(\tau_l)\vee\lambda(\tau_r)=
   \lambda(\tau'_l)\vee\lambda(\tau'_r)=\lambda(\tau')$. Similar arguments suffice when $m<k$. This completes the proof.
\quad\QED\medskip 

%%%%%%%%%%%%%%%%%%%%%%%%%%%%%%%%%%%%%%%%%%%%%%%%%%%%%%%%%%%%%%%%%%%%%%%%%%%%%%%%%
\section{Crossed product decompositions for $\SSym$ and $\YSym$}\label{S:crossed}

We observe that the surjective morphisms of Hopf algebras of
Section~\ref{S:Hopfmaps} 
\begin{equation*}%\label{E:surj}
  \raisebox{-20pt}{
   \begin{picture}(133,36)(-10,-1)

          \put(49,28){$\YSym$}
    \put(-10,0){$\QSym$}    \put(107,0){$\SSym$}
    \put(45,27){\vector(-2,-1){30}} \put(45,27){\vector(-2,-1){25}}
    \put(102.5, 2){\vector(-1,0){75}}\put(102.5, 2){\vector(-1,0){70}}
    \put(115,12){\vector(-2,1){30}}\put(115,12){\vector(-2,1){25}}
  
    \put(27,24){\scriptsize$\calL$}    \put(102,22){\scriptsize$\Lambda$}
          \put(61,5){\scriptsize$\calD$}
   \end{picture} 
}
\end{equation*}
admit splittings as coalgebras, and thus $\SSym$ is a crossed
product over $\YSym$ and $\YSym$ is a crossed product over $\QSym$.
We elucidate these structures.

Recall the poset embeddings of Section~\ref{S:combinatorics}:
\[
   \begin{picture}(113,36)(0,-1)

          \put(48,28){$\calY_n$}
    \put(0,0){$\calQ_n$}    \put(97,0){$\frakS_n$}
    \put(15,12){\epsffile{figures/duhook.eps}}
    \put(15,12){\vector(2,1){30}}
    \put(18, 2){\epsffile{figures/rhook.eps}}
    \put(18, 2){\vector(1,0){75}}
    \put(65,25){\epsffile{figures/ddhook.eps}}
    \put(65,25){\vector(2,-1){30}}
  
    \put(24,24){\scriptsize$C$}    \put(78,22){\scriptsize$\gamma$}
          \put(51,6){\scriptsize$Z$}

   \end{picture} 
\]
We use them to define linear maps as follows:
 \begin{eqnarray*}
  \calC\colon \QSym\to\YSym\,, & \setM_\alpha\ \longmapsto\ M_{C(\alpha)}\,;\\
  \Gamma \colon\YSym\to\SSym\,,& M_t\  \longmapsto\ \calM_{\gamma(t)}\,;\\
  \calZ\colon\QSym\to\SSym\,,& \setM_\alpha\ \longmapsto\ \calM_{Z(\alpha)}\,.
 \end{eqnarray*}
The following theorem is immediate from the expression for the coproduct 
on the $M$-bases of $\QSym$~\eqref{E:coproduct-QSym}, $\SSym$~\eqref{E:MRstuff1}, and $\YSym$~\eqref{E:coproduct}, and the formulas for the maps $\calD$~\eqref{E:MRstuff2}, $\Lambda$~\eqref{E:maps1} and $\calL$~\eqref{E:maps2} on these bases.

\begin{thm}
 The following is a commutative diagram of injective morphisms of coalgebras 
 which split the corresponding surjections of~\eqref{E:Des-Lambda-L}.
\[
   \begin{picture}(160,38)(0,-1)
          \put(63,26){$\YSym$}
    \put(0,0){$\QSym$}    \put(130,0){$\SSym$}
    \put(30,10){\epsffile{figures/duhook.eps}}
    \put(30,10){\vector(2,1){30}}
    \put(38, 2){\epsffile{figures/rhook.eps}}
    \put(38, 2){\vector(1,0){89}}
    \put(98,22){\epsffile{figures/ddhook.eps}}
    \put(98,22){\vector(2,-1){30}}
  
    \put(37,19){\scriptsize$\calC$}\put(112,19){\scriptsize$\Gamma$}
          \put(76,6){\scriptsize$\calZ$}
   \end{picture} 
 \]
\end{thm}

We use a theorem of Blattner, Cohen, and Montgomery~\cite{BCM86},~\cite[Ch. 7]{Mo93a}.
Suppose that $\pi\colon H\to K$ is a morphism of
Hopf algebras admitting a coalgebra splitting 
$\gamma\colon K\to H$. 
Then there is a \emph{crossed product} decomposition
\[
  H\ \cong\ A\#_{c}K
\]
where $A$, a subalgebra of $H$, is the {\em left Hopf kernel} of $\pi$:
\[
  A\ :=\ \{h\in H \mid \sum h_1\ten\pi(h_2)=h\ten 1\}
\]
and the \emph{Hopf cocycle} $c\colon K\ten K\to A$ is
 \begin{equation}\label{E:cocycle}
  c(k,k')=\sum\gamma(k_1)\gamma(k'_1)S\gamma(k_2k'_2)\,.
 \end{equation}
Note that if $\pi$ and $\gamma$ preserve
gradings, then so does the rest of this structure.

The crossed product decomposition of $\SSym$ over $\QSym$ corresponding to $(\calD,\calZ)$ was described in~\cite[Sec. 8]{AS02}.
We describe the left Hopf kernels $A$  of $\Gamma\colon\SSym\to\YSym$ 
and $B$ of $\calC\colon\YSym\to\QSym$, which are graded
with components $A_n$ and $B_n$.
Let $n>0$. Recall (Section~\ref{S:com-per-tre}) that a permutation $\tau\in\frakS_n$ admits a
unique decomposition into permutations with no global descents and a tree $t\in\calY_n$ admits a
unique decomposition into progressive trees:
\[
  \tau\ =\ \tau_1\backslash\dotsb\backslash\tau_k\qquad
   t\ =\ t_1\backslash\dotsb\backslash t_l\,.
\]
We call $\tau_k$ and $t_l$ the {\it last components} of $\tau$ and $t$.
Recall that the minimum tree $1_n=\lambda(\id_n)$ is called a  right comb.

\begin{thm} 
 A basis for  $A_n$ is the
 set $\{\calM_\tau\}$ where $\tau$ runs over all permutations of $n$ 
 whose last component is not $132$-avoiding.
 A basis for $B_n$ is the set $\{M_t\}$ where $t$ runs over all trees whose last
 component is not a right comb.
 In particular,
\[
   \dim A_n\ =\ n!-\sum_{k=0}^{n-1}k!c_{n-k-1}\qquad\textrm{and}\qquad
   \dim B_n\ =\ c_n-\sum_{k=0}^{n-1}c_k\,,
\]
 where $c_k=\#\calY_k$ is the $k$th Catalan number.
\end{thm}

\noindent{\it Proof. } By the theorem of Blattner, Cohen, and Montgomery,
 we have
\[
  \SSym\cong A\#_c\YSym \quad\textrm{and}\quad \YSym\cong B\#_c\QSym\,.
\]
 In particular $\SSym\cong A\ten\YSym$ and $\YSym\cong B\ten\QSym$
 as vector spaces.
 The Hilbert series for these graded
 algebras are therefore related by
\[
  \sum_{n\geq 0} n!t^n = \Bigl(\sum_{n\geq 0} a_nt^n\Bigr)
    \Bigl(\sum_{n\geq 0}c_nt^n\Bigr) \quad\textrm{and}\quad
  \sum_{n\geq 0} c_nt^n = \Bigl(\sum_{n\geq 0} b_nt^n\Bigr)
    \Bigl(1+\sum_{n\geq 1}2^{n-1}t^n\Bigr)\,,
\]
 where $a_n:=\dim A_n$ and $b_n:=\dim B_n$.
 Using~\eqref{E:Catalan} we deduce 
 $a_n= n!-\sum_{k=0}^{n-1}k!c_{n-k-1}$, and using
 \[ 1+\sum_{n\geq 1}2^{n-1}t^n=\frac{1}{1-\sum_{n\geq 1}t^n}\]
 we deduce  $b_n= c_n-\sum_{k=0}^{n-1}c_k$ as claimed.

 The number of permutations in $\frakS_n$ which are  $132$-avoiding and have no global descents equals the number of progressive trees in $\calY_n$, which is $c_{n-1}$. Therefore,
 $a_n$ counts the number of permutations in $\frakS_n$ whose last component is not $132$-avoiding. 
 Suppose that $\tau$ is such a permutation and $\tau=\sigma\backslash\rho$ is an arbitrary decomposition. As long as $\rho\neq\id_0$, the last component of $\rho$ is the last component of $\tau$ and hence $\rho$ is not 132-avoiding.  Thus formulas~\eqref{E:MRstuff1}  and~\eqref{E:maps1} imply that
 $(\id\ten\Lambda)\Delta(\calM_\tau)=\calM_\tau\ten 1$ and so $\calM_\tau$ 
 lies in the Hopf kernel of $\Lambda$.
 Since these elements are linearly independent, they form a basis of $A_n$ as claimed.
 
 Similarly, $b_n$ counts the number of trees in $\calY_n$ whose last component is not a comb,
 and analogous arguments using~\eqref{E:maps2} and~\eqref{E:coproduct} show that if $t$ is such a tree, then 
 $M_t$ lies in the Hopf kernel of $\calL$.
\quad\QED\medskip

%%%%%%%%%%%%%%%%%%%%%%%%%%%%%%%%%%%%%%%%%%%%%%%%%%%%%%%%%%%%%%%%%%%%%
\section{The dual of $\YSym$ and the non-commutative Connes-Kreimer Hopf algebra}\label{S:dual}

We turn now to the structure of the Loday-Ronco Hopf algebra $\LR$, which we define as the graded dual of the Hopf algebra $\YSym$.
It is known that these graded Hopf algebras are isomorphic~\cite{Foi02a,Foi02b,HNT03,Ho03,Laan}. An explicit isomorphism is obtained as the composite~\cite[Thm. 4]{HNT03},~\cite[Thm. 34]{HNT04}
\[\LR\map{\Lambda^*}(\SSym)^*\cong \SSym\map{\Lambda}\YSym\]
where the isomorphism $\SSym^*\cong \SSym$ sends an element $\calF_\sigma^*$ of the dual of the fundamental basis to $\calF_{\sigma^{-1}}$.
Thus the results of this section also apply to $\YSym$ itself.

\subsection{The dual of the monomial basis}\label{S:dualM}

Let $\{M^*_t\mid t\in\calY_\infty\}$ be the basis of $\LR$ dual to the monomial basis of $\YSym$. 

By Theorem~\ref{T:cofree}, $\YSym$ is cofree as a coalgebra and 
$\{M_t\mid t\in\calY^1\}$ is a basis for its primitive elements.
This and the form of the coproduct have the following consequence, which
also appears in~\cite[Thm. 29]{HNT04}.

\begin{thm}\label{T:free}
  $\LR$ is the free associative algebra generated by 
 $\{M^*_t\mid t\in\calY^1\}$ with
\[
  M^*_s \cdot M^*_t\ =\ M^*_{s\backslash t}\qquad 
   \mbox{for }s,t\in\calY_\infty\,.
\]
\end{thm}

The description of the coproduct requires a definition.

\begin{defi}\label{D:admissible}
 A subset $\setR$ of internal nodes of a planar binary tree is 
 {\it admissible} if for any node $x\in\setR$,  the left child $y$ of $x$ and all the descendants of $y$ are in $\setR$. Thus  any internal node  in the left
 subtree above $x$ also lies in $\setR$.
 An admissible set $\setR$ of internal nodes in a planar binary tree gives
 rise to a {\em pruning:} cut each branch connecting a node from $\setR$ to a node in 
 its complement $\setR^c$.
 For example, here is a planar binary tree whose internal nodes are labeled
 $a,b,\dotsc,h$ with an admissible set of nodes $\setR=\{f,g,d,b,c\}$.
 The corresponding pruning is indicated by the dotted line.
\begin{equation}\label{E:admissible}
\raisebox{-45pt}{
 \setlength{\unitlength}{.7142857143pt}
 \begin{picture}(175,125)
%%%%%%fig2dev -Leps -m 0.25 cutPBT.fig cutPBT.eps
   \put(0,0){\epsffile{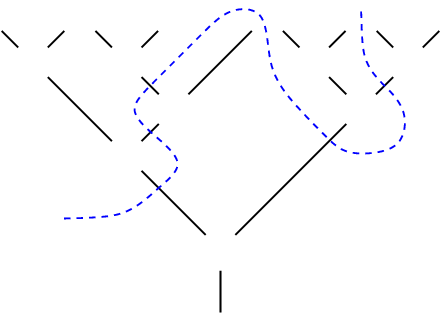}}
   \put(161, 99){$a$}
   \put(124, 99){$b$}
   \put( 46, 99){$d$}
   \put(  8, 99){$f$}
   \put(142, 81){$c$}
   \put( 66, 81){$e$}
   \put( 47, 62){$g$}
   \put( 85, 22){$h$}
 \end{picture}
}
\end{equation}
 The branches removed in such a pruning of a planar binary tree $r$ form a forest of planar binary
 trees $r_1,\dotsc,r_p$, ordered from left to right by the positions of their leaves in $r$.
 Assemble these into a planar binary tree
 $r'_\setR:=r_1\backslash r_2\backslash\dotsb\backslash r_p$.
 In~\eqref{E:admissible}, here is the pruned forest and the resulting tree $r'_\setR$:
\[
 \setlength{\unitlength}{.7142857143pt}
 \begin{picture}(170,87)(0,38)
%%%%%% fig2dev -Leps -m 0.25 cutLeft.fig cutLeft.eps
   \put(0,38){\epsffile{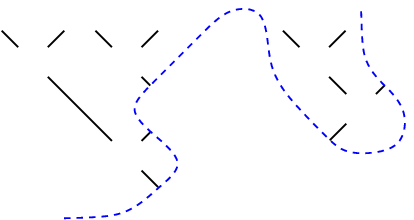}}
   \put(124, 99){$b$}
   \put( 47, 99){$d$}
   \put(  8, 99){$f$}
   \put(142, 81){$c$}
   \put( 47, 61){$g$}
 \end{picture}
\qquad
 \begin{picture}(170,92)(-40,0)
   \put(-40,40){$r'_\setR\ =$}
%%%%%% fig2dev -Leps -m 0.35 Left_tree.fig Left_tree.eps
   \put(0,0){\epsffile{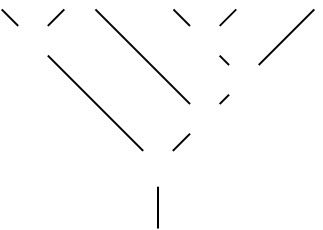}}
   \put( 80, 74){$b$}
   \put( 78, 42){$d$}
   \put(  8, 74){$f$}
   \put( 95, 57){$c$}
   \put( 59, 25){$g$}
 \end{picture}
\]
 The rest of the tree $r$ also forms a forest, which is assembled into a tree in
 a different fashion.
 If a tree $s$ in that forest is above another tree $t$ (in the original tree
 $r$) and there are no intervening components, then there is a unique leaf of $t$
 that is below the root of $s$. 
 Attach the root of $s$ to that leaf of $t$.
 As $\setR$ is admissible, there will be a unique tree in this forest below
 all the others whose root is the root of the planar binary tree
 $r''_\setR$ obtained by this assembly.
 In~\eqref{E:admissible}, here is the forest that remains and the tree
 $r''_\setR$.
\[
 \setlength{\unitlength}{.7142857143pt}
 \begin{picture}(149,125)(26,0)
%%%%%% fig2dev -Leps -m 0.35 cutRight.fig cutRight.eps
   \put(26,0){\epsffile{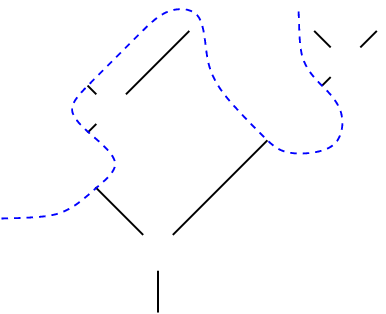}}
   \put(161, 99){$a$}
   \put( 67, 81){$e$}
   \put( 86, 22){$h$}
 \end{picture}
\qquad
 \begin{picture}(110,105)(-14,-20)
   \put(-14,30){$r''_\setR\ =$}
%%%%%% fig2dev -Leps -m 0.35 Right_tree.fig Right_tree.eps
   \put(26,0){\epsffile{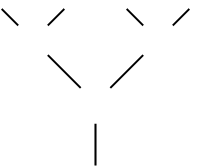}}
   \put( 87, 49){$a$}
   \put( 36, 49){$e$}
   \put( 61, 22){$h$}
 \end{picture}
\]
\end{defi}

%%%%%% This is wrong but seems irrelevant anyway
%In forming the tree $r'_\setR$, each pruned component may be further
%decomposed into its progressive components (the irreducible components for 
%the
%operation $\backslash$, as in Section~\ref{S:com-per-tre}).
%The root nodes of these progressive components determine $\setR$, and they 
%are
%distinguished as the nodes in $\setR$ that are not in the left branch above
%any other node of $\setR$.
%Call these the {\it minimal nodes in $\setR$}.
%In the example above~\eqref{E:admissible}, the minimal nodes are 
%$\{g,d,c\}$.
%%%%%% 
We record how this construction behaves under the grafting operation on trees.

\begin{lem}\label{L:attach}
 Let $s,t\in\calY_\infty$ be planar binary trees and 
 $r=s\vee  t$.
 Let $\setR$ be an admissible subset of the internal nodes of $r$, and $\setS$
 (respectively $\setT$) those nodes of $\setR$ lying in $s$ (respectively $t$).

 If the root node of $r$ lies in $\setR$, then all the nodes of $s$ lie in
 $\setR$, and 
\[ 
   r'_\setR\ =\ s \vee  t'_\setT
   \qquad\textrm{and}\qquad
   r''_\setR\ =\ t''_\setT\ .
\]
 If If the root node of $r$ does not lie in $\setR$, then 
\[ 
   r'_\setR\ =\ s'_\setS\backslash t'_\setT
   \qquad\textrm{and}\qquad
   r''_\setR\ =\ s''_\setS \vee t''_\setT\ .
\]
\end{lem}

We describe the coproduct of $\LR$ in terms of the basis
$\{M^*_t\mid t\in\calY_\infty\}$.

\begin{thm}\label{T:dual_coprod}
  For any tree $t\in\calY_\infty$, 
\begin{equation}\label{E:admissible_coprod}
  \Delta (M^*_t)\ =\ \sum M^*_{t'_\setS}\otimes M^*_{t''_\setS}\,,
\end{equation}
 the sum over all admissible subsets $\setS$ of internal nodes of $t$.
\end{thm}

We begin our proof of Theorem~\ref{T:dual_coprod}.
Recall the product formula of Theorem~\ref{T:combprod}. For
 $s\in\calY_p$ and $t\in\calY_q$,
\[
  M_s\cdot M_t\ =\ \sum_{r\in\calY_{p+q}} \alpha^r_{s,t}\,M_r\,,
\]
where $\alpha^r_{s,t}$ enumerates the set
\[
  \{\zeta\in\frakS^{(p,q)}\ \mid\ (s,t)\in\calY_p\times\calY_q
    \textrm{ is maximum such that }f_\zeta(s,t)\leq r\}
\]
and $f_\zeta\colon\calY_p\times\calY_q\to\calY_{p+q}$ is the map
\[
  f_\zeta(s,t)\ =\ \lambda\big(\gamma(s)/\gamma(t)\cdot \zeta^{-1}\big)\,.
\]
Dualizing this formula gives a formula for the coproduct of  $\LR$.
 Let $r\in\calY_n$.  Then
\begin{equation}\label{E:dual_coprod}
  \Delta(M^*_r)\ =\ \sum_{p=0}^n \, \sum_\zeta M^*_s\otimes M^*_t\,,
\end{equation}
 where $(s,t)$ is the maximum element of $\calY_p\times\calY_q$
 such that $f_\zeta(s,t)\leq r$, and the inner sum is over all 
 $\zeta\in\frakS^{(p,q)}$ such that 
 $\{(s',t')\in\calY_p\times\calY_q\mid f_\zeta(s',t')\leq r\}\neq\emptyset$.

We will deduce~\eqref{E:admissible_coprod} from~\eqref{E:dual_coprod}. 
Key to this is another reformulation of the coproduct intermediate
between these two.

For a subset $\setR\subseteq[n]$ with $p$ elements, let $\setR^c:=[n]-\setR$ be
its complement and set $q:=n-p$.
Write the elements of $\setR$ and $\setR^c$ in order
\[
  \setR\ =\ \{R_1<R_2<\dotsb<R_p\}\qquad\textrm{and}\qquad
  \setR^c\ =\ \{R_1^c<R_2^c<\dotsb<R_q^c\}\,.
\]
Define the permutation $\pi_\setR\in\frakS_n$ by 
\[
  \pi_\setR^{-1}\ :=\ (R_1,R_2,\dotsc,R_p,\,R^c_1,R^c_2,\dotsc,R^c_q)\ \in\
  \frakS^{(p,q)}\,.
\]
Then $\pi_\setR(R_i)=i$ and $\pi_\setR(R^c_i)=p+i$. Any $\zeta\in \frakS^{(p,q)}$ is of the form $\pi_\setR^{-1}$ for a unique $\setR\subseteq[n]$.

Let $\setR$ be a subset of $[n]$ as above. For a permutation $\rho\in\frakS_n$ and a tree $r\in\calY_n$ define
\[\rho|_\setR\ :=\ \st\bigl(\rho(R_1), \rho(R_2),\dotsc,\rho(R_p)\big) \qquad\text{and}\qquad r|_\setR\ :=\ \lambda\big(\gamma(r)|_\setR\bigr)\,.\]

% \[ r|_\setR\ :=\ \lambda\big(\gamma(r)(R_1), 
% \gamma(r)(R_2),\dotsc,\gamma(r)(R_p)\big)\in\calY_p\,.\]
%Here, we extend the definition of the map $\lambda\colon\frakS_n\to\calY_n$
%to give a tree $\lambda(\alpha)$ for any sequence
%$\alpha$ of distinct numbers.
%Set $\lambda(\emptyset)=|$.  Let $n\geq 1$ and assume that $\lambda$ is defined for all sequences of length less than $n$. Given a sequence $\alpha$ of
%length $n$, let $j$ be the position of the maximum element, and write $\alpha=\beta.\alpha_j.\gamma$, where the sequence $\beta$ has
%length $j-1$ and $\gamma$ has length $n-j$.
%Set
%%
% \begin{equation}\label{E:new_lambda}
%   \lambda(\alpha)\ :=\  
%     \lambda(\beta)\,\vee\,\lambda(\gamma)\,.
% \end{equation}
%%
%Note that $\lambda(\alpha)=\lambda(\st(\alpha))$.

\begin{lem}\label{L:gamma-res}
 For any $\setR\subseteq[n]$ and $r\in\calY_n$,
\begin{equation*}%\label{E:gamma-res}
\gamma(r)|_\setR=\gamma(r|_\setR)\,.
\end{equation*}
\end{lem}

\noindent{\it Proof. } Let $\sigma:=\gamma(r)|_\setR$. Since $\gamma(r)$ is
$132$-avoiding (Section~\ref{S:combinatorics}), so is $\sigma$. Hence 
$\sigma=\gamma\bigl(\lambda(\sigma)\bigr)$, and 
$\gamma(r)|_\setR=\gamma\Bigl(\lambda\bigl(\gamma(r)|_\setR\bigr)\Bigr)=\gamma(r|_\setR)$.
\quad\QED\medskip

\begin{thm}\label{T:intermediate}
 Let $r\in\calY_n$. 
 Then
\begin{equation}\label{E:intermediate}
   \Delta(M^*_r)\ =\ \sum_{\substack{\setR\subseteq[n]\\\lambda(\pi_\setR)\leq r}}
      M^*_{r|_\setR}\otimes M^*_{r|_{\setR^c}}\ .
\end{equation}
\end{thm}

\noindent{\it Proof. }
 Let $\setR\subseteq[n]$ with $\#\setR=p$ and $\zeta:=\pi_\setR^{-1}\in\frakS^{(p,q)}$.
 Since the map $f_\zeta$ is order-preserving, the minimum element in its image
 is $f_\zeta(1_p,1_q)=\lambda(\zeta^{-1})$, and so the 
 sums in~\eqref{E:dual_coprod} and~\eqref{E:intermediate} are over the same sets.
 We only need show that if $\lambda(\zeta^{-1})\leq r$ then $(r|_\setR, r|_{\setR^c})$ is maximum among those pairs
 $(s,t)\in\calY_p\times\calY_q$ such that $f_\zeta(s,t)\leq r$.
  We first establish the corresponding fact about permutations; namely that
 if $\zeta^{-1}\leq\rho$ then $(\rho|_\setR, \rho|_{\setR^c})$ is maximum among those pairs
 $(\sigma,\tau)\in\frakS_p\times\frakS_q$ such that $\varphi_\zeta(\sigma,\tau)\leq \rho$.

 The permutation
 $\upsilon:=\varphi_\zeta(\sigma,\tau)=(\sigma/\tau)\cdot\zeta^{-1}$ satisfies
\[
   \upsilon(R_i)\ =\ \sigma(i)\qquad\textrm{and}\qquad
   \upsilon(R^c_j)\ =\ p+\tau(j)\,,
\]
 for $i=1,\dotsc,p$ and $j=1,\dotsc,q$.
 Thus $\upsilon|_\setR=\sigma$ and $\upsilon|_{\setR^c}=\tau$.
  We describe the inversion set of $\upsilon$:
 \begin{eqnarray*}
   (R_i,R_j)\in\Inv(\upsilon)    &\Longleftrightarrow& (i,j)\in\Inv(\sigma)\\
   (R^c_i,R^c_j)\in\Inv(\upsilon)&\Longleftrightarrow& (i,j)\in\Inv(\tau)\\
   (R^c_i,R_j)\in\Inv(\upsilon)  &\Longleftrightarrow& R^c_i<R_j
 \end{eqnarray*}
 There are no inversions of $\upsilon$ of the form $(R_i,R^c_j)$.

The above includes a description of $\Inv(\zeta^{-1})$ (choosing $\sigma=\id_p$, $\tau=\id_q$).
Since the weak order on $\frakS_n$ is given by inclusion of inversion sets,
 we see that  for a permutation $\rho\in\frakS_n$,
\[
  \zeta^{-1}\leq\rho\ \Longleftrightarrow\ 
  \{(R^c_i,R_j)\mid R^c_i<R_j\}\ \subseteq\ \Inv(\rho)\,.
\]
 Since $(i,j)$ is an inversion of $\rho|_\setR$ is and only if 
 $(R_i,R_j)$ is an inversion of $\rho$, we see that if $\zeta^{-1}\leq\rho$,
 then $(\rho|_\setR,\rho|_{\setR^c})$ is maximum among all pairs 
 $(\sigma,\tau)\in\frakS_p\times\frakS_q$ such that
 $(\sigma/\tau)\cdot\zeta^{-1}\leq\rho$. 

We finish the proof by deducing the fact about trees. Suppose $\lambda(\zeta^{-1})\leq r$. Let $\rho:=\gamma(r)$. Then
$\zeta^{-1}\leq\rho$, by Theorem~\ref{T:galois}. Suppose $f_\zeta(s,t)\leq r$.
Let $\sigma:=\gamma(s)$ and $\tau=\gamma(t)$. Then 
$\lambda\bigl((\sigma/\tau)\cdot\zeta^{-1}\bigr)=f_\zeta(s,t)\leq r$, so
$\varphi_\zeta(\sigma,\tau)\leq\rho$. By the fact about permutations,
$\sigma\leq\rho|_\setR$ and $\tau\leq\rho|_{\setR^c}$. Applying $\lambda$ we obtain $s\leq r|_\setR$ and $t\leq r|_{\setR^c}$. It remains to verify that
$f_\zeta(r|_\setR,r|_{\setR^c})\leq r$. This is equivalent to 
$\varphi_\zeta\bigl(\gamma(r|_\setR),\gamma(r|_{\setR^c})\bigr)\leq \rho$.
This holds by the fact about permutations, since $\gamma(r|_\setR)=\rho|_\setR$ and $\gamma(r|_{\setR^c})=\rho|_{\setR^c}$
by Lemma~\ref{L:gamma-res}. 
\quad\QED\medskip

\begin{lem}\label{L:graft-res}
 Let  $1\leq j\leq n$, $\sigma\in\frakS_{j-1}$, $\tau\in\frakS_{n-j}$,
  $s\in\calY_{j-1}$, and $t\in\calY_{n-j}$.  Let $\setR\subseteq[n]$, $\setS:=\setR\cap[1,j-1]$, and $\setT:=\setR\cap[j+1,n]$. Then
 \begin{equation*}%\label{E:graft-res}
(\sigma\vee\tau)|_\setR\ =\ \begin{cases}
(\sigma|_\setS)\vee(\tau|_{\setT-j}) & \text{ if $j\in\setR$,}\\
(\sigma|_\setS)\backslash(\tau|_{\setT-j}) & \text{ if $j\notin\setR$;}
\end{cases}
\qquad
(s\vee t)|_\setR\ =\ \begin{cases}
(s|_\setS)\vee (t|_{\setT-j}) & \text{ if $j\in\setR$,}\\
(s|_\setS)\backslash (t|_{\setT-j}) & \text{ if $j\notin\setR$.}
\end{cases}
\end{equation*}
\end{lem}

\noindent{\it Proof. } The statement for permutations is immediate from the definitions. Applying $\gamma$ to both sides of the remaining equality, and using~\eqref{E:def-gamma},~\eqref{E:gamma-slash}, and Lemma~\ref{L:gamma-res}
we deduce the statement for trees.
\quad\QED\medskip

We complete the proof of Theorem~\ref{T:dual_coprod} by
showing that under a natural labeling of the internal nodes
of a tree $r\in\calY_n$, admissible subsets of nodes are exactly those subsets
$\setR\subseteq[n]$ such that $\lambda(\pi_\setR)\leq r$, and that given such a subset
$\setR$, 
\[
   r|_\setR\ =\ r'_\setR\qquad\textrm{and}\qquad
   r|_{\setR^c}\ =\ r''_\setR\,.
\]

Label the $n$ internal nodes  of a tree $r\in\calY_n$ with the integers
$1,2,\dotsc,n$ in the following recursive manner. Write $r=s\vee t$ with $s\in\calY_{j-1}$ and $t\in\calY_{n-j}$, $1\leq j\leq n$. 
Assume the nodes of $s$ and $t$ have been labeled. The root node of $r$ is labeled with $j$, if a node comes from $s$, it retains its
label, and if a node in comes from $t$, we increase its label by $j$.
Note that the label of any internal node of $r$ is bigger than
all the labels of nodes in its left subtree and smaller  than all the labels of nodes in its right subtree.
%With this labeling, the sums of Theorems~\ref{T:dual_coprod}
%and~\ref{T:intermediate} coincide. 

\begin{lem}
 Let $\setR\subseteq[n]$ and $r\in\calY_n$.
 We consider $\setR$ to be a subset  of internal nodes of $r$, under the above
 labeling. 
 Then
\[
  \lambda(\pi_\setR)\  \leq\  r\quad\Longleftrightarrow\quad
  \setR\textrm{ is admissible.}
\]
\end{lem}

\noindent{\it Proof. }
 Let $\setR\subseteq[n]$ and $r\in\calY_n$, and set $\rho:=\gamma(r)$.
By Theorem~\ref{T:galois}, $\lambda(\pi_\setR)\leq r\iff \pi_\setR\leq\rho$.
 In the proof of Theorem~\ref{T:intermediate}, we showed that 
\[
  \pi_\setR\ \leq\ \rho\ 
   \Longleftrightarrow\ \textrm{whenever }i<j\textrm{ with }
    i\not\in\setR,\textrm{ and }j\in\setR, 
    \textrm{ then }\rho(i)>\rho(j)\,.
\]
 Equivalently, if $i<j$ with $\rho(i)<\rho(j)$ and $j\in\setR$, then 
 $i\in\setR$.
To show that this is equivalent to $\setR$ being admissible, we only need to verify
  that if $i<j$ with $\rho(i)<\rho(j)$, then
 in $r$ the node labeled $i$ is in the left subtree above the node labeled $j$.  

Let $h$ be the label of the root node of $r$, $1\leq h\leq n$. Thus $r=s\vee t$ with $s\in\calY_{h-1}$ and $t\in\calY_{n-h}$.
 By definition of $\gamma$~\eqref{E:def-gamma}, $\rho=\gamma(r)=\gamma(s)\vee\gamma(t)$. By definition of grafting of permutations~\eqref{E:def-grafting-perm}, $\rho(h)=n$.  Thus $\rho$ achieves its maximum on the label $h$ of the root node. Suppose $j=h$. By construction of the labeling, all labels  $i<j$ belong to the left subtree above the root,
 which shows that the claim holds in this case. If $j\neq h$, since $\rho$ is $132$-avoiding, we must have either $i<j<h$ or $h<i<j$.  In the former case,
 $\rho(i)<\rho(j)\iff\gamma(s)(i)<\gamma(s)(j)$; in the latter,
 $\rho(i)<\rho(j)\iff\gamma(t)(i-h)<\gamma(t)(j-h)$.
  The claim now follows by induction on $n$.
\quad\QED\medskip

 The following lemma completes the proof of Theorem~\ref{T:dual_coprod}.

\begin{lem}
 Let $\setR\subseteq[n]$ be an admissible subset of nodes of a tree
 $r\in\calY_n$, labeled as above.
 Then
\[
   r|_\setR\ =\ r'_\setR\qquad\textrm{and}\qquad
   r|_{\setR^c}\ =\ r''_\setR\,.
\]
\end{lem}

\noindent{\it Proof. }
Write  $r=s\vee t$ with $s\in\calY_{j-1}$ and $t\in\calY_{n-j}$.
Thus $j$ is  the label of the root node of $r$, the set of labels of the nodes of $s$ and $t$ are respectively $[1,j{-}1]$ and $[j{+}1,n]$.

%Let $\setS:=\setR\cap[1,j{-}1]$ and $\setT:=\setR\cap[j{+}1,n]$.

 Suppose that $j\in\setR$.
 As $\setR$ is admissible, $[1,j{-}1]\subseteq\setR$, and by Lemma~\ref{L:graft-res},
\[
   r|_\setR\ =\ s \vee (t|_{\setR\cap[j{+}1,n]-j})
    \qquad\textrm{and}\qquad
   r|_{\setR^c}\ =\ t|_{\setR^c\cap[j{+}1,n]-j}\,.
\]
Proceeding by induction we may assume that $t|_{\setR\cap[j{+}1,n]-j}=t'_{\setR\cap[j{+}1,n]-j}$ and $t|_{\setR^c\cap[j{+}1,n]-j}=t''_{\setR\cap[j{+}1,n]-j}$.
Together with Lemma~\ref{L:attach} this gives $r|_\setR=r'_\setR$ and
   $r|_{\setR^c}=r''_\setR$.

 Similarly, if $j\not\in\setR$, then by Lemma~\ref{L:graft-res},
\[
   r|_\setR\ =\ (s|_{\setR\cap[1,j{-}1]}) \backslash  (t|_{\setR\cap[j{+}1,n]-j})
    \qquad\textrm{and}\qquad
   r|_{\setR^c}\ =\ (s|_{\setR^c\cap[1,j{-}1]}) \vee (t|_{\setR^c\cap[j{+}1,n]-j})\,.
\]
Induction and an application of Lemma~\ref{L:attach} complete the proof.
\quad\QED\medskip

\subsection{$\LR$ and the non-commutative Connes-Kreimer Hopf algebra}\label{S:LR_NCK}

We use Theorems~\ref{T:free} and~\ref{T:dual_coprod} to give an explicit isomorphism between
$\LR$ and the non-commutative Connes-Kreimer Hopf algebra, $\NCK$ of
Foissy~\cite[Sec. 5]{Foi02a}.  Holtkamp constructed a less explicit isomorphism~\cite[Thm. 2.10]{Ho03}. Palacios~\cite[Sec. 4.4.1]{Pa02} obtained an explicit description of this isomorphism which is equivalent to ours. Foissy~\cite{Foi02b} showed that the two Hopf algebras are isomorphic by exhibiting a dendriform structure on $\NCK$.

As an algebra, $\NCK$ is freely generated by the set
of all finite rooted planar trees. 
Monomials of rooted planar trees are naturally identified with {\em ordered forests} (sequences of rooted planar trees), so  $\NCK$ has a linear basis of such forests. The identity element corresponds to the empty forest $\emptyset$.
The algebra $\NCK$ graded by the total number of nodes in a forest.
Here are some forests.
\[
%%%%%%fig2dev -Leps -m 0.4 f1.fig f1.eps  
  \epsffile{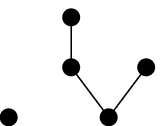}\ ,\qquad\quad
  \epsffile{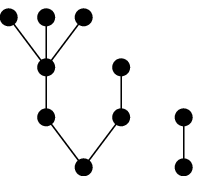}\ ,\qquad\quad
  \epsffile{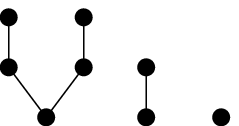}\ .
\]
An subset $\setR$ of nodes of a forest is {\it admissible}  if for any node $x\in\setR$, every node above $x$ also lies in
$\setR$. 
Given an admissible subset of nodes in a forest $f$, we prune the
forest by removing the edges connecting nodes of $\setR$ to nodes of its
complement.
The pruned pieces give a planar forest $f'_\setR$, and the pieces that remain
also form a forest, $f''_\setR$.
For example, here is a pruning of the third forest above, and the resulting forests:
\begin{equation}\label{E:cut}
%fig2dev -Leps -m 0.4 fCut.fig fCut.eps
\raisebox{-18pt}{
  \epsffile{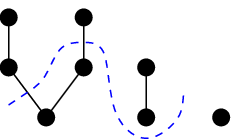}\qquad
  \raisebox{15pt}{$f'_\setS\ =$  }\ \epsffile{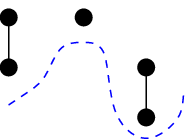}\  
  \raisebox{15pt}{$=$}\  \raisebox{5pt}{\epsffile{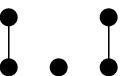}}\qquad
  \raisebox{15pt}{$f''_\setS\ =$ }\ 
  \raisebox{5pt}{\epsffile{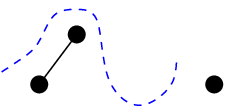}}\ 
  \raisebox{15pt}{$=$ }\ \raisebox{5pt}{\epsffile{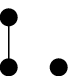}}\ \raisebox{6pt}{.}}
\end{equation}

The coproduct in $\NCK$ is given by
\begin{equation}\label{E:NCK_coprod}
  \Delta (f)\ =\ \sum f'_\setS\ten f''_\setS\,,
\end{equation}
the sum over all admissible subsets $\setS$ of nodes of the forest $f$.
To prove $\NCK\cong \LR$, we furnish a bijection $\varphi$
between planar forests $f$ of rooted planar trees with $n$ nodes and planar
binary trees with $n$ internal nodes that preserves these structures.

We construct $\varphi$ recursively. Set $\varphi(\emptyset):=\,|\,$.
Removing the root from a planar rooted tree $t$ gives a planar forest $f$,
and we set $\varphi(t):=\varphi(f)/\epsffile{figures/1.eps}$. Finally, 
given a forest $f=(t_1,t_2,\dotsc,t_n)$, where each $t_i$ is a planar rooted
tree, set
$\varphi(f):=\varphi(t_1)\backslash\varphi(t_2)\backslash\dotsb\backslash\varphi(t_n)$.

For example, $\varphi({\scriptstyle \bullet})=\epsffile{figures/1.eps}$,
$\varphi(\barBell)\ =\ \varphi({\scriptstyle \bullet})/\epsffile{figures/1.eps}\ 
 =\ \epsffile{figures/12.eps}$, and so
 \begin{eqnarray*}
   \varphi\left(\raisebox{-9pt}{\epsfysize=20pt\epsffile{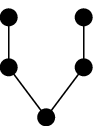}}\right)
    &=&
   \big(\epsffile{figures/12.eps}\,\backslash\,\epsffile{figures/12.eps}\big)
    \,/\,\epsffile{figures/1.eps} 
   \ =\ \raisebox{-5pt}{\epsffile{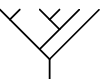}}
    \\
   \varphi\left(\raisebox{-9pt}{\epsfysize=20pt\epsffile{figures/f3.eps}}\right)
   &=& \raisebox{-5pt}{\epsffile{figures/34125.eps}}\,
   \backslash\,\epsffile{figures/12.eps}
   \,\backslash\,\epsffile{figures/1.eps}\ =\ 
  \raisebox{-5pt}{\epsffile{figures/67458231.eps}}\,.
\end{eqnarray*}

The last example above shows that  the planar binary tree of~\eqref{E:admissible} and the forest of~\eqref{E:cut} correspond to each other
under $\varphi$. Compare the admissible subsets and prunings illustrated in~~\eqref{E:admissible} and~\eqref{E:cut}.
Under $\varphi$ the images of the nodes above a node $x$
consist of all the internal nodes in the left branch above the image of $x$. 
Thus admissible subsets of nodes of a forest $f$ correspond to admissible
subsets of internal nodes of the planar binary tree $\varphi(f)$.
Similarly, the assembly of the pieces given by a cut corresponding to admissible
sets also correspond, as may be seen from these examples and Lemma~\ref{L:attach}. 

We deduce an isomorphism between the 
 non-commutative Connes-Kreimer Hopf algebra and the Loday-Ronco Hopf algebra.

Define a linear map $\Phi:\NCK\to\LR$ by
\begin{equation}\label{E:def-Phi}
\Phi(f)\ :=\ M^*_{\varphi(f)}
\end{equation}

\begin{thm}\label{T:NCK-LR}
 The map $\Phi$ 
 is an isomorphism of graded Hopf algebras $\NCK\cong\LR$.
 \end{thm}
\noindent{\it Proof. }
 Theorem~\ref{T:free} guarantees that $\Phi$ is a morphism of algebras and the
 preceding discussion shows that
 $\Phi$ is a morphism of coalgebras. It is easy to see that $\varphi$ is a bijection between the set of ordered forests with $n$ nodes and the set of
 planar binary trees with $n$ internal nodes. Thus $\Phi$ is an isomorphism of graded Hopf algebras.
\quad\QED\medskip

The non-commutative Connes-Kreimer Hopf algebra  carries a canonical involution.
Given a plane forest $f$, let $f^{\rm r}$ be its reflection across a
vertical line on the plane. It is clear that
\[ (f^{\rm r})^{\rm r}=f\,, \qquad
  (f\cdot g)^{\rm r}\ =\ g^{\rm r}\cdot f^{\rm r}\,,
  \qquad\mbox{\rm and}\qquad
  \Delta (f)^{{\rm r}\otimes{\rm r}}\ =\ \Delta (f^{\rm r})\,;
\]
in other words, the map $f\mapsto f^{\rm r}$ is an involution, an algebra anti-isomorphism, and a coalgebra isomorphism of the noncommutative Connes-Kreimer  Hopf algebra with itself.

We deduce the existence of a canonical involution on $\YSym$, which we construct recursively. Define $|^{\rm r}:=|$. For a progressive tree $t$, write
$t=s\vee |$, and define $t^{\rm r}:=s^{\rm r}\vee |$. Finally, for an arbitrary
 planar binary tree $t$, consider its decomposition into progressive trees $t=t_1\backslash t_2 \backslash \cdots \backslash t_k$ (Section~\ref{S:com-per-tre}) and define
\[
   t^{\rm r}\ :=\  (t_k)^{\rm r}\backslash \cdots \backslash
   (t_2)^{\rm r} \backslash  (t_1)^{\rm r}\,.
\]
For instance,
\[
 %fig2dev -Leps -m0.2  10786934125.fig 10786934125.eps
  \mbox{\epsffile{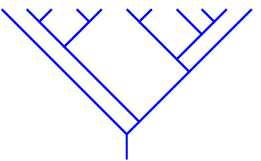}}^{\rm r}
   \qquad\raisebox{15pt}{$=$}\qquad
  %fig2dev -Leps -m0.2  98671042351.fig 98671042351.eps
  \epsffile{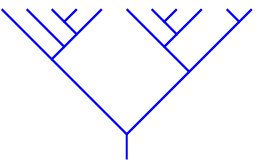}
\]

\begin{cor} The map $\YSym\to\YSym$, $M_t\mapsto M_{t^{\rm r}}$, is an involution, an algebra isomorphism, and a coalgebra anti-isomorphism.
 \end{cor}
\noindent{\it Proof. } By construction, $\varphi(f^{\rm r})=\varphi(f)^{\rm r}$. We may thus transport the result from $\NCK$ to $\LR$ via $\Phi$ (and to $\YSym$ via duality).
 \quad\QED\medskip

Since the map  $t\mapsto t^{\rm r}$ does not preserve the Tamari order on $\calY_n$, the involution does not admit a simple expression on the $\{F_t\}$-basis of $\YSym$. We also remark that there is a commutative diagram

\[\begin{picture}(80,55)(0,0)
\put(0,50){$\YSym$}  \put(70,50){$\YSym$}
    \put(0,0){$\QSym$}    \put(70,0){$\QSym$}
    \put(15,45){\vector(0,-1){33}}
    \put(85,45){\vector(0,-1){33}}
    \put(35, 2){\vector(1,0){33}}
    \put(35, 52){\vector(1,0){33}}
     \put(48,4){\rm r}    \put(48,54){\rm r}
      \put(0,25){$\calD$} \put(90,25){$\calD$}
  \end{picture} \]

The bottom map sends $M_\alpha\mapsto M_{\alpha^{\rm r}}$, where $\alpha^{\rm
  r}=(a_k,\ldots,a_2,a_1)$ is the reversal of the composition
$\alpha=(a_1,a_2,\ldots,a_k)$.  This is an involution, an algebra isomorphism,
and a coalgebra anti-isomorphism of $\QSym$ with itself. The map $\alpha\mapsto
\alpha^{\rm r}$ {\em is} order-preserving, so the involution on $\QSym$ is also
given by $F_\alpha\mapsto F_{\alpha^{\rm r}}$. 

\subsection{Symmetric functions and the Connes-Kreimer Hopf algebra}\label{S:Sym-CK}

Let $\CK$ be the Connes-Kreimer Hopf algebra. 
It is the free commutative algebra
generated by the set of all finite rooted (non-planar) trees. Commutative
monomials of rooted trees are naturally identified with {\em unordered forests}
(multisets of rooted trees), so $\CK$ has a linear basis consisting of such
unordered forests. 
The coproduct of $\CK$ is defined in terms of admissible subsets of nodes in the
same way as for $\NCK$~\eqref{E:NCK_coprod}. 
$\CK$ is a commutative graded Hopf algebra.

Given an ordered forest $f$ of planar trees, let $U(f)$ be the unordered forest
obtained by forgetting the left-to-right order among the trees in $f$, and the
left-to-right ordering among the branches emanating from each node in each tree
in $f$.  
The map $U:\NCK\to\CK$ is a surjective morphism of Hopf algebras.

Consider rooted trees  in which each node has at most one child.
These  are sometimes called {\em ladders}. Let $\ell_n$ be the ladder with $n$ nodes. Clearly,
\[\Delta(\ell_n)=\sum_{i=0}^n \ell_i\otimes\ell_{n-i}\,.\]
It follows that the subalgebra of $\CK$ generated by $\{\ell_n\}_{n\geq 0}$, is a Hopf
subalgebra, isomorphic to the Hopf algebra of symmetric functions via the map
\[\Sym\inc\CK\,, \qquad h_n\mapsto \ell_n\,.\]
Here $h_n$ denotes the {\em complete} symmetric function. 

Recall that the graded dual of the Hopf algebra of quasi-symmetric functions is the Hopf algebra of non-commutative symmetric functions: $\QSym^*=\NSym$.
Dualizing the map $\calL:\YSym\onto\QSym$ (Proposition~\ref{P:ontoHopf})
we obtain an injective morphism of
Hopf algebras, which by~\eqref{E:maps2} is given by
\[\NSym\inc\YSym\,,\qquad M_\alpha^*\mapsto M_{C(\alpha)}^*\,.\]
By definition of the map $C$~\eqref{E:def-C} and Theorem~\ref{T:free}, if $\alpha=(a_1,\ldots,a_k)$ then
\[M_{C(\alpha)}^*=M_{1_{a_1}}^*\cdots M_{1_{a_k}}^*\,.\]
The bijection $\varphi$ of Section~\ref{S:LR_NCK} sends the ladder $\ell_n$ (viewed as a planar rooted tree) to the comb $1_n$. Therefore, composing with the isomorphism of Theorem~\ref{T:NCK-LR} we obtain an injective morphism of Hopf algebras
\[\NSym\inc\NCK\,,\qquad M_\alpha^*\mapsto \ell_{a_1}\cdots\ell_{a_k}\,.\]

The canonical map $\NSym\onto\Sym$ sends $M_\alpha^*$ to the complete
symmetric function $h_{a_1}\cdots h_{a_k}$. We have shown:

\begin{thm} There is a commutative diagram of graded Hopf algebras
\begin{equation}\label{E:CK-Sym}
 \raisebox{-20pt}{  \begin{picture}(80,60)(0,0)
\put(-3,50){$\NSym$}  \put(70,50){$\NCK$}
    \put(2,0){$\Sym$}    \put(73,0){$\CK$}
    \put(15,45){\vector(0,-1){33}} \put(15,45){\vector(0,-1){29}}
    \put(85,45){\vector(0,-1){33}}  \put(85,45){\vector(0,-1){29}}
    \put(35, 2){\vector(1,0){33}}\put(35, 2){\epsffile{figures/rhook.eps}}
    \put(35, 52){\vector(1,0){33}}\put(35, 52){\epsffile{figures/rhook.eps}}
      \end{picture} }
\end{equation}
  \end{thm}

Let $\GL:=\CK^*$ denote the graded dual of the Connes-Kreimer Hopf algebra.
As shown by Hoffman~\cite{Ho02}, this is the (cocommutative) Hopf algebra of rooted trees  constructed by Grossman and Larson in~\cite{GL89}.
Dualizing~\eqref{E:CK-Sym} we obtain the following commutative diagram of graded Hopf algebras:
\begin{equation}\label{E:CK-Sym-dual}
 \raisebox{-20pt}{ \begin{picture}(80,55)(0,0)
\put(5,50){$\GL$}  \put(70,50){$\YSym$}
    \put(2,0){$\Sym$}    \put(70,0){$\QSym$}
    \put(15,45){\vector(0,-1){33}}\put(15,45){\vector(0,-1){29}}
    \put(85,45){\vector(0,-1){33}}\put(85,45){\vector(0,-1){29}}
    \put(30, 2){\vector(1,0){38}}\put(30, 2){\epsffile{figures/rhook.eps}}
    \put(30, 52){\vector(1,0){38}}\put(30, 52){\epsffile{figures/rhook.eps}}
      \end{picture} }
\end{equation}

\def\cprime{$'$}
\providecommand{\bysame}{\leavevmode\hbox to3em{\hrulefill}\thinspace}
\providecommand{\MR}{\relax\ifhmode\unskip\space\fi MR }
\providecommand{\MRhref}[2]{%
  \href{http://www.ams.org/mathscinet-getitem?mr=#1}{#2}
}
\providecommand{\href}[2]{#2}

\end{document}